\documentclass[11pt]{article}
\usepackage{amssymb,amsmath,mathrsfs,amsfonts,graphicx,epsf}
\usepackage{color}
\usepackage{xcolor}
\usepackage{lpic}

\DeclareMathAlphabet\gothic{U}{euf}{m}{n}
 \input xy
 \xyoption{all}

\usepackage{fullpage}

\newcommand{\immagine}[3][16]{
\begin{figure}[htb]
\begin{center}
\includegraphics[width=#1cm]{#2.pdf}
\caption{#3}
\label{fig:#2}
\end{center}
\end{figure}}
\renewcommand{\r}[1]{(\ref{#1})}

\newcommand{\srarc}{sR-arclength}

\newcommand{\ul}[1]{{\mathbf #1}}
\renewcommand{\th}{\theta}

\usepackage{mathrsfs}
\usepackage{graphicx}
\usepackage{amssymb}
\usepackage{color}
\newtheorem{theorem}{Theorem}
\newtheorem{corollary}[theorem]{Corollary}
\newtheorem{lemma}[theorem]{Lemma}
\newtheorem{proposition}[theorem]{Proposition}
\newtheorem{definition}[theorem]{Definition}
\newtheorem{remark}[theorem]{Remark}
\newenvironment{proof}{{\it Proof.}~~}{\hfill$\square$}
\newcommand{\bt}{\begin{theorem}}
\newcommand{\et}{\end{theorem}}
\newcommand{\bl}{\begin{lemma}}
\newcommand{\el}{\end{lemma}}
\newcommand{\bp}{\begin{proposition}}
\newcommand{\ep}{\end{proposition}}
\newcommand{\bc}{\begin{corollary}}
\newcommand{\ec}{\end{corollary}}
\newcommand{\bdeff}{\begin{definition}}
\newcommand{\edeff}{\end{definition}}
\newcommand{\brem}{\begin{remark}}
\newcommand{\erem}{\end{remark}}
\newcommand{\bproof}{\begin{proof}}
\newcommand{\eproof}{\end{proof}}

\newcommand{\bi}{\begin{itemize}}
\newcommand{\iii}{\item}
\renewcommand{\i}{\item}
\newcommand{\ei}{\end{itemize}}
\newcommand{\bd}{\begin{description}}
\newcommand{\ed}{\end{description}}
\newcommand{\be}{\begin{enumerate}}
\newcommand{\ee}{\end{enumerate}}

\newcommand{\bqn}{\begin{eqnarray}}
\newcommand{\eqn}{\end{eqnarray}}
\newcommand{\eqnn}{\nonumber\end{eqnarray}}
\newcommand{\eqnl}[1]{\label{#1}\end{eqnarray}}

\newcommand{\nn}{\nonumber}
\newcommand{\ba}[1]{\begin{array}{#1}}
\newcommand{\ea}{\end{array}}

\newcommand{\R}{\mathbb{R}}

\newcommand{\lam}{\lambda}
\newcommand{\g}{q(.)}
\newcommand{\eps}{\varepsilon}

\def\cn{\mathrm{cn\,}}
\def\sn{\mathrm{sn\,}}
\def\dn{\mathrm{dn\,}}
\def\E{\mathrm{E\,}}
\def\Z{{\mathbf{Z}}}

\newcommand{\sgn}{\textrm{sgn}\,}
\newcommand{\tcut}{T_{cut}}
\newcommand{\tconj}{T_{conj}}
\newcommand{\tcusp}{T_{cusp}}
\newcommand{\tpend}{T_{pend}}

\newcommand{\EXP}{\textsf{Exp}}
\DeclareMathOperator{\Cut}{Cut}
\DeclareMathOperator{\Con}{Con}

\newcommand{\pmec}{${\bf (P_{MEC})}$}
\newcommand{\pcurve}{${\bf(P_{curve})}$}
\newcommand{\pprojective}{${\bf (P_{projective})}$}
\newcommand{\pproj}{\pprojective}
\newcommand{\VP}{${\bf (VP)}$}

\newcommand{\PTR}{PT\R^2}
\newcommand{\Pt}[1]{\left( #1 \right)}
\newcommand{\Pg}[1]{\left\{ #1 \right\}}
\newcommand{\Pq}[1]{\left[ #1 \right] }
\newcommand{\x}{{\bf x}}
\renewcommand{\b}[1]{{\bf #1}}
\renewcommand{\kappa}{K}

\title{Curve cuspless reconstruction via sub-Riemannian geometry}
\author{U. Boscain, R. Duits, F. Rossi, Y. Sachkov}

\begin{document}

\maketitle

\begin{abstract}
We consider the problem of minimizing $\int_{0}^\ell \sqrt{\xi^2 +K^2(s)}\, ds $
for a planar curve having fixed initial
and final positions and directions. The total length $\ell$ is free. Here $s$ is the variable of arclength parametrization, $K(s)$ is the
curvature of the curve and $\xi>0$ a parameter. This
problem comes from a model of geometry of vision due to
Petitot, Citti and Sarti.

We study existence of local and global minimizers for this problem. We prove that if for a certain choice of boundary conditions there is no global minimizer, then there is neither a local minimizer nor a geodesic.

We finally give properties of the set of boundary conditions for which there exists a solution to the problem.
\end{abstract}

\section{Introduction}
In this paper we are interested in the following variational problem\footnote{In this paper, by $S^1$ we mean $\R/\sim$ where $\theta\sim\theta'$ if $\theta=\theta'+2n\pi$, $n\in\Z$. By $P^1$ we mean $\R/\approx$ where $\theta\approx\theta'$ if $\theta=\theta'+n\pi$, $n\in\Z$.}:
\bd
\iii{\bf (P)}
Fix $\xi>0$ and  $(x_{in},y_{in},\theta_{in}),(x_{fin},y_{fin},\theta_{fin})\in \R^2\times S^1$.
On the space of (regular enough)
planar curves, parameterized by plane-arclength\footnote{Here by plane-arclength we mean the arclength in $\R^2$, for which we use the variable $s$. Later on, we consider also parameterizations by arclength on $\R^2\times S^1$ or $\R^2\times P^1$, that we call sub-Riemannian arclength (\srarc\ for short), for which we use the variable $t$. We will also use the variable $\tau$ for a general parametrization.}
find the solutions of:
\bqn
&&\ul{x}(0)=(x_{in},y_{in}),~~~\ul{x}(\ell)=(x_{fin},y_{fin}) ,\nn\\
&&\dot{\ul{x}}(0)=(\cos(\theta_{in}), \sin(\theta_{in})) ,~~\dot{\ul{x}}(\ell)=(\cos(\theta_{fin}), \sin(\theta_{fin})),\nn\\
&&    \int_0^\ell\sqrt{\xi^2+\kappa^2(s)}~ds\to\min~~~~~ (\mbox{here $\ell$ is free.})
\eqn
\ed
Here $\kappa=\frac{\dot x \ddot  y-\dot  y\ddot  x}{(\dot  x^2+\dot  y^2)^{3/2}}$ is the geodesic curvature of the planar  curve $\ul{x}(\cdot)=(x(\cdot),y(\cdot))$. This problem comes from a model proposed by Petitot, Citti and Sarti (see \cite{citti-sarti,petitot,petitot-libro, sanguinetti} and references therein) for the mechanism of reconstruction of corrupted curves used by the visual cortex V1. The model is explained in detail in Section \ref{s-petitot}.

It is convenient to formulate the problem {\bf (P)} as a problem of optimal control, for which the functional spaces are also more naturally specified.
\bd
\iii{\pcurve}
Fix $\xi>0$ and  $(x_{in},y_{in},\theta_{in}),(x_{fin},y_{fin},\theta_{fin})\in \R^2\times S^1$.
In the space of integrable controls $v(\cdot):[0,\ell]\to\R$,
find the solutions of:
\bqn
&&(x(0),y(0),\theta(0))=(x_{in},y_{in},\theta_{in}),~~~(x(\ell),y(\ell),\theta(\ell))=(x_{fin},y_{fin},\theta_{fin}) ,\nn\\
&&\Pt{\ba{c}
\frac{dx}{ds}(s)\\
\frac{dy}{ds}(s)\\
\frac{d\theta}{ds}(s)
\ea
}=\Pt{\ba{c}
\cos(\th(s))\\
\sin(\th(s))\\
0
\ea}+v(s)\Pt{\ba{c}
0\\
0\\
1
\ea}\label{e-pcurve}\\
&&    \int_0^\ell\sqrt{\xi^2+\kappa(s)^2}~ds= \int_0^\ell\sqrt{\xi^2+v(s)^2}~ds\to\min ~~~(\mbox{here $\ell$ is free})
\eqn
\ed
Since in this problem we are taking $v(\cdot)\in L^1([0,\ell])$, we have that the curve $\g=(x(\cdot),y(\cdot),\th(\cdot)):[0,\ell]\to\R^2\times S^1$ is absolutely continuous and the planar curve  $\ul{x}(\cdot):=(x(\cdot),y(\cdot)):[0,\ell]\to\R^2$ is in $W^{2,1}([0,\ell])$.

\brem
Notice that the function $\sqrt{\xi^2+K^2}$ has the same asymptotic behaviour, for $K\to0$ and for $K\to\infty$ of the function
$\phi(K)$ introduced by Mumford and Nitzberg in their functional for image segmentation (see \cite{mumford21}).
\erem

The main issues we address in this paper are related to existence of minimizers for problem \pcurve. More precisely, for \pcurve\ the first question we are interested in is the following:
\bd
\iii[Q1)] Is it true that for every initial and final condition, the problem \pcurve\ admits a {\bf global minimum}?
\ed
In \cite{suzdal}  it was shown that there are initial and final conditions for which \pcurve\  does not admit a minimizer. More precisely, it was shown that there exists a minimizing sequence for which the limit is a trajectory not satisfying the boundary conditions. See Figure \ref{fig:angoli-inizio}.

\immagine[10]{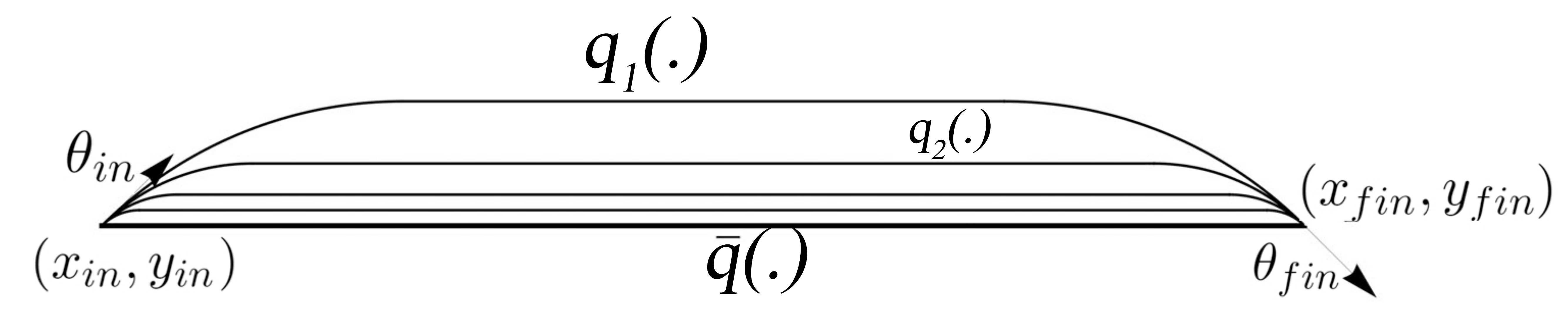}{Minimizing sequence $q_n$ converging to a non-admissible curve $\bar q$ (angles at the beginning/end).}

From the modeling point of view, the non-existence of global minimizers is not a crucial issue. It is very natural to assume that the visual cortex looks only for local minimizers, since it seems reasonable to expect that it primarly compares nearby trajectories. Hence, a second problem we address in this paper is the existence of local minimizers 
for the problem \pcurve. More precisely, we answer the following question:
\bd
\iii[Q2)] Is it true that for every initial and final condition the problem \pcurve\ admits a {\bf local minimum}? If not, what is the set of boundary conditions for which a local minimizer exists?
\ed

\medskip
The main result of this paper is the following.
\bt Fix an initial and a final condition $q_{in}=(x_{in},y_{in},\th_{in})$ and $q_{fin}=(x_{fin},y_{fin},\th_{fin})$ in $\R^2\times S^1$. The only two following cases are possible.
\begin{enumerate}
\iii  There exists a solution (global minimizer) for \pcurve\ from $q_{in}$ to $q_{fin}$.
\iii The problem  \pcurve\ from $q_{in}$ to $q_{fin}$ does not admit neither a global nor a local minimum nor a geodesic.
\end{enumerate}
Both cases occur, depending on the boundary conditions.
\label{t-maini}
\et

We recall that a curve $\g$ is a geodesic if for every sufficiently small interval $[t_1,t_2]\subset Dom(\g)$, the curve $\g_|{_{[t_1,t_2]}}$ is a minimizer between $q(t_1)$ and $q(t_2)$.

One of the main interests of \pcurve\ is that it admits minimizers that are in $W^{1,1}$ but are not Lipschitz, as we will show in Section \ref{ch:lav}. As a consequence, controls lie in $L^1$ but not in $L^\infty$. This is an interesting phenomenon for control theory: indeed, to find minimizers, one usually applies the Pontryagin Maximum Principle (PMP in the following), that is a generalization of the Euler-Lagrange condition. But the standard formulation of the PMP holds for $L^\infty$ controls; this obliges us to use a generalization of the PMP for \pcurve, that we discuss in Section \ref{s-genpmp}. Details of this interesting aspect of \pcurve\ are given in Section \ref{ch:lav}. This also explains the reason for which we need to define variational problems, global and local minimizers in the space $W^{1,1}$, see Section \ref{s-sR}.

The second sentence of {\bf Q2} is interesting, since one could compare the limit boundary conditions for which a mathematical reconstruction occurs with the limit boundary conditions for which a reconstruction in human perception experiments is observed. Indeed, it is well known from human perception experiments that the visual cortex V1 does not connect all initial and final conditions, see e.g. \cite{petitot-libro}. With this goal, we have computed numerically the configurations for which a solution exists, see Figure~\ref{fig:Q2}.
\immagine{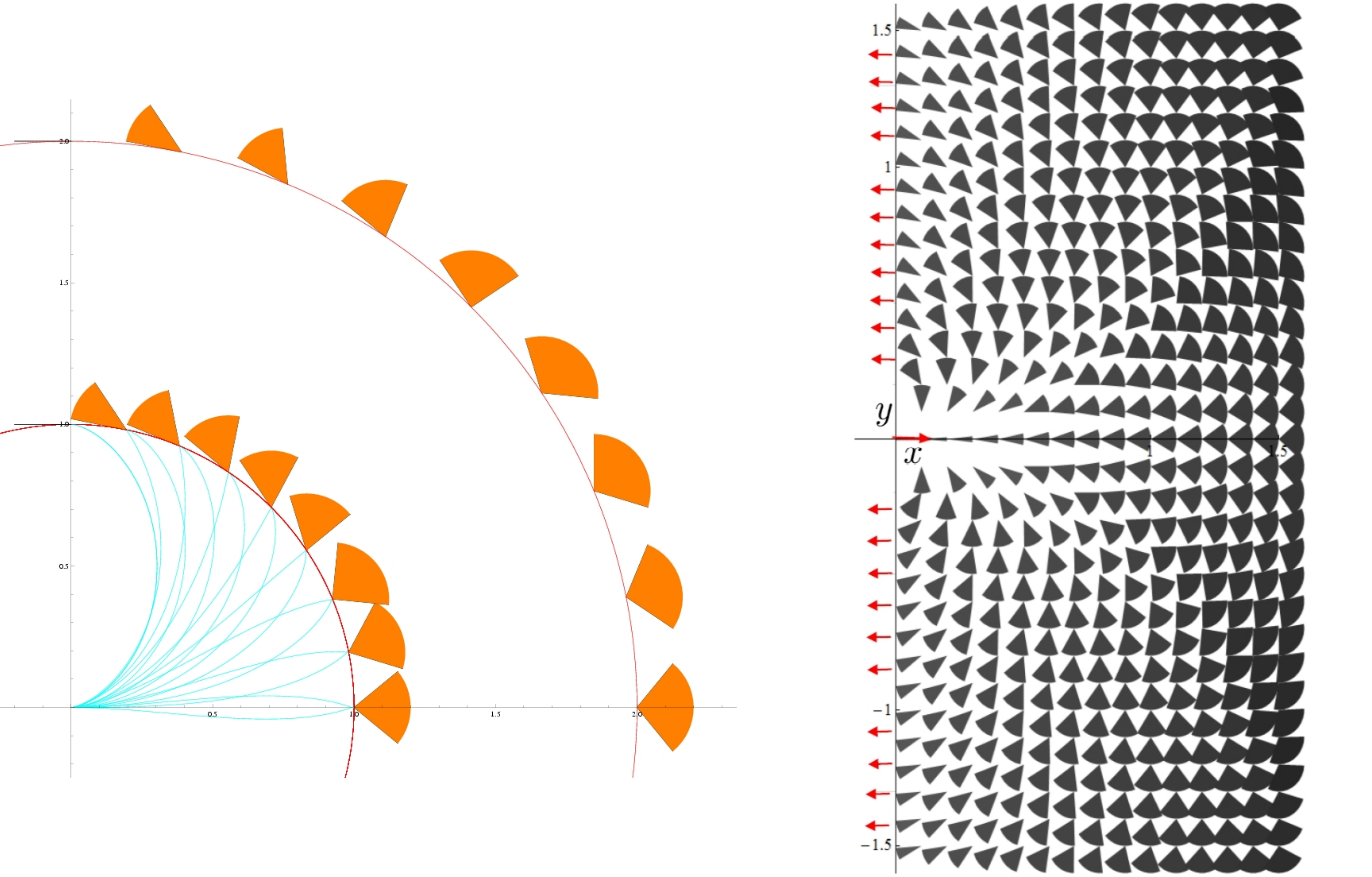}{Configurations for which we have existence of minimizers with $\xi=1$. For other $\xi\neq 1$, one can recover the corresponding figure via dilations, as explained in Remark \ref{remark:xi}. Due to invariance of the problem under rototranslations on the plane, one can always assume that $q_{in}=(0,0,0)$. Left: We study the cases $x^2_{fin}+y^2_{fin}=1$ and $x^2_{fin}+y^2_{fin}=4$, with $y\geq0$. The case $y\leq 0$ can be recovered by symmetry. In the case $x^2_{fin}+y^2_{fin}=1$ minimizing curves are also shown. Right: For each point on the right half-plane, we give the set of configurations for which we have existence of minimizers.}

The structure of the paper is as follows. In Section \ref{s-petitot} we briefly describe the model by Petitot-Citti-Sarti for the visual cortex V1. We state it as a problem of optimal control (more precisely a sub-Riemannian problem), that we denote by \pproj. The problem \pcurve\ is indeed a modified version of \pproj. In Section \ref{s-sR} we recall definitions and main results in sub-Riemannian geometry, that is the main tool we use to prove our results. In Section \ref{s-pmec} we define an auxiliary mechanical problem (crucial for our study), that we denote with \pmec, and study the structure of geodesics for it. In Section \ref{s-problems} we describe in detail the relations between problems \pcurve, \pprojective\ and \pmec, with an emphasis on the connections between the minimizers of such problems. In Section \ref{s-main} we prove the main results of the paper, i.e. Theorem \ref{t-maini}.

\section{The model by Petitot-Citti-Sarti for V1}
\label{s-petitot}

In this section, we recall a model describing how the human visual cortex V1 reconstructs curves which are partially hidden or corrupted. The goal is to explain the connection between reconstruction of curves and the problem \pcurve\ studied in this paper.

The model we present here was initially due to Petitot \cite{petitot,petitot-libro}, based on previous work by Hubel-Wiesel \cite{hubel} and Hoffman \cite{hoffman}, then refined by Citti and Sarti \cite{citti-sarti, sanguinetti}, and by the authors of the present paper in \cite{nostro-vision,duits-q1,duits-q2}. It was also studied by Hladky and Pauls in \cite{hladky}.

In a simplified model\footnote{For example, in this model we do not take into account the fact that the continuous space of stimuli is implemented via a discrete set of neurons.} (see \cite[p. 79]{petitot-libro}), neurons of V1 are grouped into {\it orientation columns}, each of them being sensitive to visual stimuli at a given point  of the retina and for a given direction  on it. The retina is modeled by the real plane, i.e. each point is represented by $(x,y)\in\R^2$, while the directions at a given point are modeled by the projective line, i.e. $\theta\in P^1$. Hence, the primary visual cortex V1 is modeled by the so called {\it projective tangent bundle} $\PTR:=\R^2\times P^1$. From a neurological point of view, orientation columns are in turn grouped into {\it hypercolumns}, each of them being sensitive to stimuli at a given point $(x,y)$ with any direction. In the same hypercolumn, relative to a point $(x,y)$ of the plane, we also find neurons that are sensitive to other stimuli properties, like colors, displacement directions, etc...  In this paper, we  focus only on directions and therefore  each hypercolumn is represented by a fiber $P^1$ of the bundle $\PTR$.
Orientation columns are connected between them in two different ways. The first kind is given by  {\it vertical connections}, which connect
 orientation columns belonging to the same hypercolumn and sensible to similar directions. The second is given by the  {\it horizontal connections}, which connect orientation columns in different (but not too far) hypercolumns and sensible to the same directions. See Figure \ref{fig:f-hyper-bis}.

\immagine{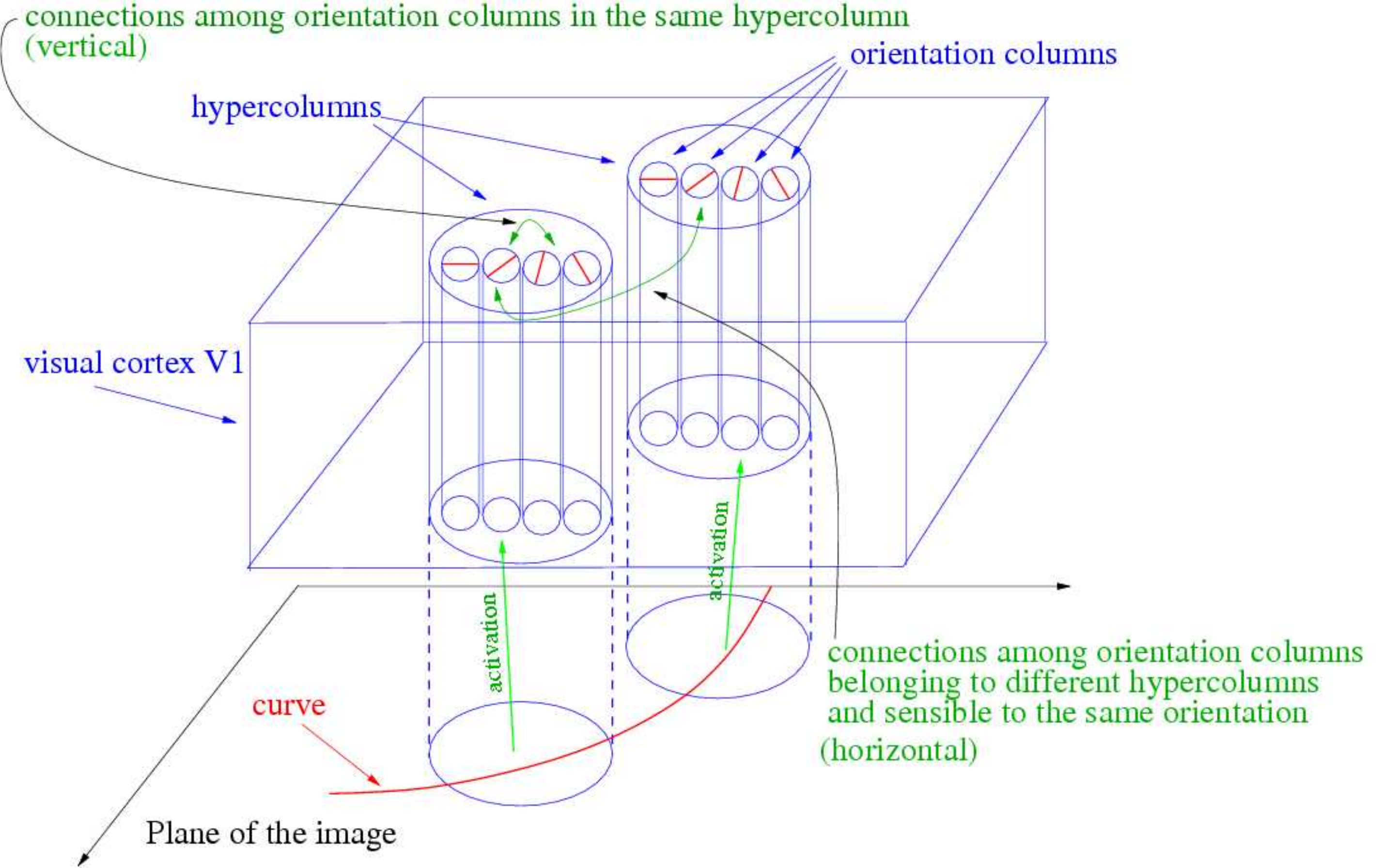}{A scheme of the primary visual cortex V1.}

In other words, when V1 detects a (regular enough) planar curve $(x(\cdot),y(\cdot)):[0,T]\to\R^2$ it computes a ``lifted curve'' in  $\PTR$ by including a new variable $\theta(\cdot):[0,T]\to P^1$, defined in $W^{1,1}([0,T])$, which  satisfies:
\bqn
&&\Pt{\ba{c}
\frac{d x}{d\tau}(\tau)\\
\frac{d y}{d\tau}(\tau)\\
\frac{d \th}{d\tau}(\tau)
\ea
}=u(\tau)\Pt{\ba{c}
\cos(\th(\tau))\\
\sin(\th(\tau))\\
0
\ea}+v(\tau)\Pt{\ba{c}
0\\
0\\
1
\ea} \mbox{~~~for some ~~~}u,v:[0,T]\rightarrow\R.
\label{eq-contrSR}
\eqn
The new variable $\theta(.)$ plays the role of the direction in $P^1$ of the tangent vector to the curve. Here it is natural to require  $u(\cdot),v(\cdot)\in L^1([0,T])$.  This specifies also which regularity we need for the planar curve to be able to compute its lift: we need a curve in $W^{2,1}$.

Observe that the lift is not unique in general: for example, in the case in which there exists an interval $\Pq{\tau_1,\tau_2}$ such that $\frac{d x}{d\tau}(\tau)=\frac{d y}{d\tau}(\tau)=0$ for all $\tau\in\Pq{\tau_1,\tau_2}$, one has to choose $u=0$ on the interval, while the choice of $v$ is not unique. Nevertheless, the lift is unique (modulo $L^1$) in many relevant cases, e.g. if $\frac{d x}{d\tau}(\tau)^2+\frac{d y}{d\tau}(\tau)^2= 0$ for a finite number of times $\tau\in[0,T]$.

\b{In the following we call a planar curve a {\it liftable curve} if it is in $W^{2,1}$ and its lift is unique.}

Consider now a  liftable curve $(x(\cdot),y(\cdot)):[0,T]\to\R^2$ which is interrupted in an interval $(a,b)\subset[0,T]$. Let us call
$(x_{in},y_{in}):=(x(a),y(a))$ and $(x_{fin},y_{fin}):=(x(b),y(b))$. 
Notice that the limits $\theta_{in}:=\lim_{\tau\uparrow a}\theta(\tau)$ and  $\theta_{fin}:=\lim_{\tau\downarrow b}\theta(\tau)$ are well defined, since $\th(.)$ is an absolutely continuous curve. In the model by Petitot, Citti, Sarti and the authors of the present article \cite{nostro-vision,citti-sarti,pinwheel}, the visual cortex reconstructs the curve by minimizing the energy necessary to activate orientation columns which are not activated by the curve itself. This is modeled by the minimization of the functional
\bqn
&&J=\int_a^b \left(\xi^2 u(\tau)^2+v(\tau)^2\right)~d\tau\to\min \label{eq-minQ1},~~~\mbox{(here $a$ and $b$ are fixed)}\label{e-5}.
\eqn
Indeed, $\xi^2 u(\tau)^2$ (resp. $v(\tau)^2$) represents the (infinitesimal) energy necessary to activate horizontal (resp. vertical) connections. The parameter
$\xi>0$ is used to fix the relative weight of the horizontal and vertical connections, which have different phisical dimensions.  The minimum is  taken on the set of curves which are solution of (\ref{eq-contrSR}) for some $u(\cdot),v(\cdot)\in L^1([a,b])$ and satisfying boundary conditions
\bqn
(x(a),y(a),\theta(a))=(x_{in},y_{in},\theta_{in}),~~~(x(b),y(b),\theta(b))=(x_{fin},y_{fin},\theta_{fin}).
\eqnn

Minimization of the cost (\ref{e-5}) is equivalent to the minimization of the cost (which is invariant by reparameterization)
\bqn
{\cal L}=\int_a^b \sqrt{\xi^2 u(\tau)^2+v(\tau)^2}~d\tau=\int_a^b \|\dot{\ul{x}}(\tau)\| \sqrt{\xi^2+\kappa(\tau)^2}~d\tau, \label{e-7}
\eqnn 
where $\ul{x}=(x,y)$ and with $b>a$ fixed. See a proof of such equivalence in \cite{yuri1}.
%

We thus define the following problem:
\bd
\iii{\pprojective:}
Fix $\xi>0$ and  $(x_{in},y_{in},\theta_{in}),(x_{fin},y_{fin},\theta_{fin})\in \R^2\times P^1$. In the space of integrable controls $u(\cdot),v(\cdot):[0,T]\to\R$, find the solutions of:
\bqn
&&(x(0),y(0),\theta(0))=(x_{in},y_{in},\theta_{in}),~~~(x(T),y(T),\theta(T))=(x_{fin},y_{fin},\theta_{fin}) ,\nn\\
&&\Pt{\ba{c}
\frac{d x}{d\tau}(\tau)\\
\frac{d y}{d\tau}(\tau)\\
\frac{d \th}{d\tau}(\tau)
\ea
}=u(\tau)\Pt{\ba{c}
\cos(\th(\tau))\\
\sin(\th(\tau))\\
0
\ea}+v(\tau)\Pt{\ba{c}
0\\
0\\
1
\ea}\nn\\
&&    {\cal L}= \int_0^T\sqrt{\xi^2u(\tau)^2+v(\tau)^2}~d\tau=\int_0^T\|\dot{\ul{x}}(\tau)\| \sqrt{\xi^2+\kappa(\tau)^2}~d\tau \to\min ~~~(\mbox{here $T$ is free})
\eqnn
\ed
Observe that here $\theta\in P^1$, i.e. angles are considered without orientation\footnote{Notice that the vector field $(\cos\th,\sin\th,0)$ is not continuous on $\PTR$. Indeed, a correct definition of \pproj\ needs two charts, as explained in detail in \cite[Remark 12]{nostro-vision}. In this paper, the use of two charts is implicit, since it plays no crucial role.}.

The optimal control problem \pprojective\ is well defined. Moreover, it is a sub-Riemmanian problem, see Section \ref{s-sR}. We have remarked in \cite{nostro-vision} that a solution always exists. We have also studied a similar problem in \cite{nostro-PTS2}, when we deal with curves on the sphere $S^2$ rather than on the plane $\R^2$.

One the main interests of \pprojective\ is the possibility of associating to it a hypoelliptic diffusion equation which can be used to reconstruct images (and not just curves),  and for contour enhancement. This point of view was developed in \cite{nostro-vision,citti-sarti,duits-q1,duits-q2}.

However, its main drawback (at least for the problem of reconstruction of curves) is the existence of minimizers with cusps, see e.g. \cite{suzdal}. 
Roughly speaking, cusps are singular points in which velocity changes its sign. More formally, we say that a curve trajectory $(q(.),(u(.),v(.)))$ has a cusp at $\bar \tau\in[0,T]$ if  $u(\tau)$ changes its sign in a neighbourhood\footnote{More precisely, it exists $\eps>0$ such that $u(a)u(b)<0$ for almost every $a\in \Pt{\bar \tau-\eps,\bar \tau}, b\in \Pt{\bar \tau,\bar\tau+\eps}$.} of $\bar \tau$. Notice that in a neighborhood of a cusp point, the tangent direction (with no orientation) is well defined. A minimizer with cusps is represented in Figure \ref{fig:f-cusp}.

\immagine[8]{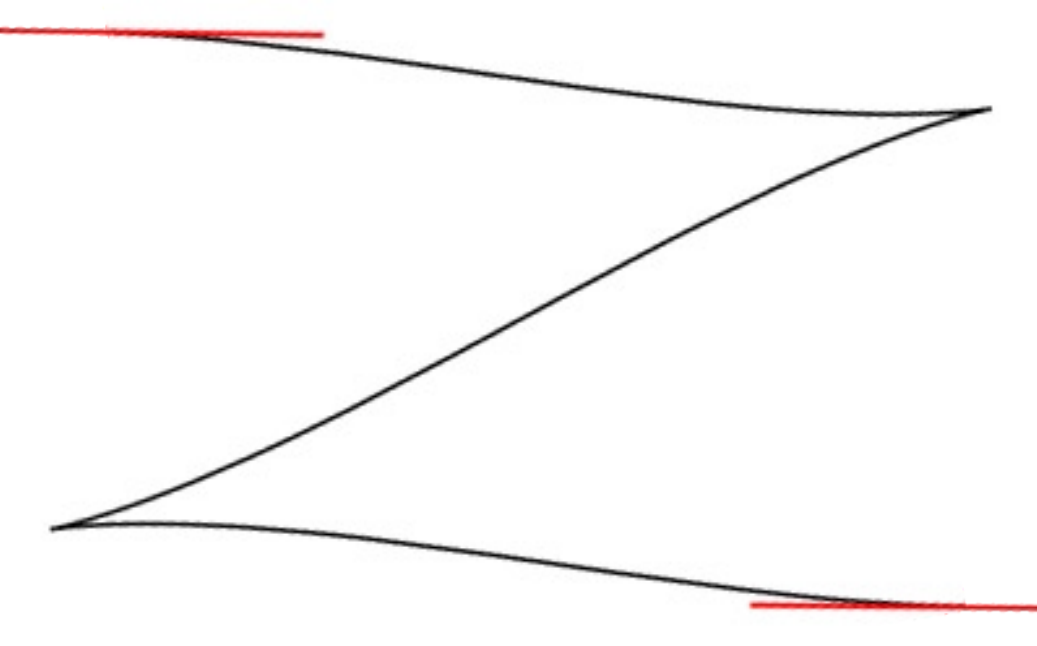}{A minimizer with two cusps.}

However, to our knowledge, the presence of cusps has not been observed in human perception experiments, see e.g. \cite{petitot-libro}. For this reason, people started looking for a way to require that no trajectories with cusps appear as solutions of the minimization problem. In \cite{citti-sarti,duits-q2} the authors proposed to require trajectories parameterizated by spatial arclength, i.e. to impose $\|\dot{\ul{x}}\|=u=1$. In this way cusps cannot appear. Notice that assuming $u=1$,  directions must be considered with orientation, since now  the direction of $\dot{\ul{x}}$ is defined in $S^1$. In fact, cusps are precisely the points where the spatial arclength parameterization breaks down. By fixing $u=1$ we get the optimal control problem \pcurve\ on which this paper is focused. 

\brem \label{r-cuspdistribuita} We also define an ``angular cusp'' as follows: we say that a pair trajectory-control $(q(.),(u(.),v(.)))$ has an ``angular cusp'' at $\bar \tau\in[a,b]$  if there exist a neighbourhood $B:=(\bar \tau-\eps,\bar \tau+\eps)$ such that $u(\tau)\equiv 0$ on $B$ and $\th(\bar \tau-\eps)\neq \th(\bar \tau+\eps)$. Angular cusps are of the kind $q(\tau)=(x_0,y_0,\th_0+\int_0^\tau v(\sigma)\,d\sigma)$.

The minimum of the distance between $(x_0,y_0,\th_0)$ and $(x_0,y_0,\th_1)$ with arbitrary $\th_0,\th_1$ is realized by such kind of trajectories. This is the only interesting case in which we need to deal with such trajectories. Indeed, even assuming that a solution $\g$ of \pprojective\ satisfies $u_1\equiv 0$ on a neighbourhood  $\bar t$ only, then analyticity of the solution\footnote{Analyticity of the solution is proved below, see Remark \ref{analytic}.} implies that $u_1\equiv 0$ on the whole $\Pq{0,T}$, and hence $\g$ steers $(x_0,y_0,\th_0)$ to some $(x_0,y_0,\th_1)$.
\erem

\section{Optimal control}\label{s-sR}

In this section, we give the fundamental definitions and results from optimal control, and the particular cases of sub-Riemannian problems, that we will use in the following. For more details about sub-Riemannian geometry, see e.g. \cite{bellaiche,gromov,mont}.

\subsection{Minimizers, local minimizers, geodesics}

In this section, we give main definitions of optimal control. Observe that we deal with curves in the space $W^{1,1}$ to deal with the problem \pcurve, see Section \ref{ch:lav}.

\bdeff \label{d-vp} Let $M$ be an $n$ dimensional smooth manifold and $f_u:q\mapsto f_u(q)\in T_qM$ be a 1-parameter family of smooth vector fields depending on the parameter $u\in\R^m$. Let $f^0:M\times \R^m\rightarrow [0,+\infty)$ be a smooth function of its arguments. We call \b{variational problem} (denoted by {\VP} for short) the following optimal control problem
\bqn
&&\dot q(\tau)=f_{u(\tau)}(q(\tau)), \label{q1}\\
&&\int_{0}^T f^0(q(\tau),u(\tau))~d\tau\to\min,~~ T\mbox{~free}\label{q2}\\
&&q(0)=q_0,~~~q(T)=q_1,\label{q3}\\
&&u(\cdot)\in \cup_{T>0}L^1([0,T], \R^m),~~q(\cdot)\in \cup_{T>0}W^{1,1}([0,T],M)\label{q4}
\eqn
\edeff

Following \cite[Ch. 8]{vinter}, we endow $\cup_{T>0}W^{1,1}([0,T],M)$ with a topology.
\bdeff
Let $q_1(.),q_2(.)\in\cup_{T>0}W^{1,1}([0,T],M)$, with $q_1$ defined on $[0,T_1]$ and $q_2$ on $[0,T_2]$. Extend $q_1$ on the whole time-interval $[0,\max\Pg{T_1,T_2}]$ by defining $q_1(t):=q_1(T_1)$ for $t>T_1$, and similarly for $q_2$. We define the distance between $q_1(.)$ and $q_2(.)$ as
$$\|q_1(.)-q_2(.)\|_{W^{1,1}}:= |q_1(0)-q_2(0)|+\|\dot q_1(.)-\dot q_2(.)\|_{L^1}.$$
\edeff
From now on, we endow $\cup_{T>0}W^{1,1}([0,T],M)$ with the topology induced by this distance. It is clear that this distance is induced by the norm in $W^{1,1}$. For more details, see \cite[Ch. 8]{vinter}.

We now give definitions of minimizers for \VP.
\bdeff
\label{d-minimizers}
We say that a pair trajectory-control $(q(\cdot),u(\cdot))$ is a minimizer if it is a solution of {\VP}.\\
We say that it is a local minimizer if there exists an open neighborhood $B_{q(\cdot)}$  of $q(\cdot)$ in $\cup_{T>0}W^{1,1}([0,T],\R^m)$, endowed with the topology defined above, such that all $(\bar q(\cdot),\bar u(\cdot))$ satisfying (\ref{q1})-(\ref{q3}), with $\bar q(\cdot)\in  B_{q(\cdot)}$, have a larger or equal cost.\\
We say  that it is a geodesic if for every sufficiently small interval $[t_1,t_2]\subset Dom(q(\cdot))$, the pair  $(q(\cdot),u(\cdot))|_{[t_1,t_2]}$ is a minimizer of
$\int_{t_1}^{T} f^0(q(\tau),u(\tau))~d\tau$ from $q(t_1)$ to $q(t_2)$ with $T$ free.
\edeff
\brem It is interesting to observe that, in general, one can have the same trajectory $q(.)$ realized by two different controls $u_1(.),u_2(.)$. For this reason, one has to specifiy the control to have the cost of a trajectory. Nevertheless, for the problems studied in this article, it is easy to prove that, for a given trajectory $q(.)\in W^{1,1}([0,T],\R^m)$ satisfying \r{q1} for some control $u(.)$, then such control is unique.
\erem

In this paper we are interested in studying problems that are particular cases of {\bf (VP)}, see Section \ref{s-problems}. In particular, we study the problem \pmec\ defined in Section \ref{s-pmec}, that is a 3D contact problem (see the definition below). For such problem we apply a standard tool of optimal control, namely the Pontryagin Maximum Principle (described in the next section), and then derive properties for \pcurve\ from the solution of \pmec.

\subsection{Sub-Riemannian manifolds}
In this section, we recall the definition of sub-Riemannian manifolds and some properties of the corresponding Carnot-Caratheodory distance. We recall that sub-Riemannian problems are special cases of optimal control problems.

\renewcommand{\Delta}{\blacktriangle}

\bdeff
A sub-Riemannian manifold is a triple $(M,\Delta,{\mathbf g}),$
where
\bi
\i $M$ is a connected smooth manifold of dimension $n$;
\i $\Delta$ is a Lie bracket generating smooth distribution of constant rank $m<n$; i.e., $\Delta$
is a smooth map that associates to $q\in M$  an $m$-dim subspace $\Delta(q)$ of $T_qM$, and $\forall~q\in M$, we have
\bqn\label{Hor}
\mathrm{span}\Pg{[X_1,[\ldots[X_{k-1},X_k]\ldots]](q)~|~X_i\in\mathrm{Vec}(M)\mbox{~and~}X_i(p)\in\Delta(p)~\forall~p\in
M}=T_qM.\nonumber \eqn Here ${\rm Vec(M)}$ denotes the set of
smooth vector fields on $M$. \i ${\mathbf g}_q$ is a Riemannian
metric on $\Delta(q)$, that is, smooth as a function of $q$.
\ei
\edeff

The Lie bracket generating condition \r{Hor} is also known as the
H\"ormander\break condition, see \cite{hor}.

\bdeff \label{d-distanza}
A Lipschitz continuous curve $q(.):[0,T]\to M$ is said to be {\it horizontal} if
$\dot q(\tau)\in\Delta(q(\tau))$ for almost every $\tau\in[0,T]$.
Given a horizontal curve $q(.):[0,T]\to M$, the {\it length of $q(.)$} is
\bqn
l(q(.))=\int_0^T \sqrt{ {\mathbf g}_{q(\tau)} (\dot q(\tau),\dot q(\tau))}~d\tau.
\eqn
The {\it distance} induced by the sub-Riemannian structure on $M$ is the
function
\bqn
d(q_0,q_1)=\inf \Pg{l(q(.))\mid q(0)=q_0,\,q(T)=q_1,\, q(.)\ \mbox{\rm horizontal Lipschitz continuous curve}}.
\label{e-distanza}
\eqn
\edeff
Notice that the length of a curve is invariant by time-reparametrization of the curve itself.

The hypothesis of connectedness of $M$ and the Lie bracket generating assumption
for the distribution guarantee the finiteness and the continuity of
$d(\cdot,\cdot)$ with respect to the topology of $M$ (Rashevsky-Chow's theorem; see, for instance, \cite{agra-book}).

The function $d(\cdot,\cdot)$ is called the Carnot--Caratheodory distance. It gives to $M$ the structure of a metric space (see \cite{bellaiche,gromov}).

Observe that \pprojective\ and \pmec\ defined in Section \ref{s-pmec} are both sub-Riemannian problems. Indeed, defining $$X_{1}=\left(\ba{c}
\cos\th\\
 \sin\th\\
0
\ea\right),
~~
X_{2}=\left(\ba{c}
0\\
0\\
1
\ea\right),$$ one has that \pprojective\ is sub-Riemannian with $M=\PTR$, $\Delta_q=\mathrm{span}\Pg{X_1(q),X_2(q)}$ and ${\mathbf g(q)}$ such that $X_1(q),X_2(q)$ is an orthonormal basis. For \pmec, simply replace $M=SE(2)$.

\subsection{The Pontryagin Maximum Principle on 3D contact manifolds}
\label{s-pmp}

In the following, we state some classical results from geometric control theory which hold for the 3D contact case. For simplicity of notation, we only consider structures defined globally by a pair of vector fields, that are sometimes called {\it ``trivialized structures''}.

\bdeff[3D contact problem] Let $M$ be a 3D manifold and let $X_{1},X_{2}$ be two smooth vector fields such that dim(Span$\{X_{1},X_{2},[X_{1},X_{2}] \}(q)$)=3 for every $q\in M$. The variational problem
\bqn
\dot q=u_1X_{1}+u_2X_{2},~~q(0)=q_0,~~~q(T)=q_1,~~\int_0^T\sqrt{(u_1(\tau))^2+(u_2(\tau))^2 } d\tau\to \min
\eqnl{e:3Dprob}
is called a 3D contact problem.
\edeff
Observe that a 3D contact manifold is a particular case of a sub-Riemannian manifold, with $\Delta=\mathrm{span}\Pg{X_1,X_2}$ and ${\mathbf g}_{q(t)}$ is uniquely determined by the condition ${\mathbf g}_{q(\tau)}(X_i,X_j)=\delta_{ij}$. In particular, each 3D contact manifold is a metric space when endowed with the Carnot-Caratheodory distance.

When the manifold is analytic and the orthonormal frame can be assigned through $m$ analytic vector fields, we say that the sub-Riemannian manifold is {\bf analytic}. This is the case of the problems studied in this article.
\brem \label{r-Linf}
In the problem above the final time $T$ can be free or fixed since the cost is invariant by time reparameterization. As a consequence the spaces $L^1$ and $W^{1,1}$ in \r{q4} can be replaced with $L^\infty$ and $Lip$ (like in \r{e-distanza}), since we can always reparameterize trajectories in such a way that $u_1(\tau)^2+u_2(\tau)^2=1$ for every $\tau\in[0,T]$. If $u_1(\tau)^2+u_2(\tau)^2=1$ for a.e. $\tau\in[0,T]$  we say that the curve is parameterized by \srarc. See \cite[Section 2.1.1]{suzdal} for more details.
\erem
We now state first-order necessary conditions for our problem.
\bp[Pontryagin Maximum Principle for 3D contact problems] \label{p-pmp}
In the 3D contact case, a curve parameterized by \srarc\  is a geodesic if and only if it is the projection of a solution of the Hamiltonian system corresponding to the Hamiltonian
 \bqn
\label{eq-HH}
H(q,p)=\frac12( \langle p, X_{1}(q)\rangle^2 + \langle p, X_{2}(q)\rangle^2),~~q\in M,~p\in T^\ast_qM,
\eqn
lying on the level set  $H=1/2$.
\ep
 This simple  form of the Pontryagin Maximum Principle follows from the absence of nontrivial abnormal extremals in 3D contact geometry, as a consequence of the condition $\mathrm{dim(Span}\{X_{1},X_{2},[X_{1},X_{2}] \}(q))=3$ for every $q\in M$, see \cite{agra-exp}. For a general form of the Pontryagin Maximum Principle, see \cite{agra-book}.

\brem\label{analytic} As a consequence of Proposition \ref{p-pmp}, for 3D contact problems, geodesics and the corresponding controls are always smooth and even analytic if $M, X_{1},X_{2}$ are analytic, as it is the case for the problems studied in this article. Analyticity of geodesics in sub-Riemannian geometry holds for general analytic sub-Riemannian manifolds having no abnormal extremals. For more details about abnormal extremals, see e.g. \cite{davidenote,mont}.
\erem 
 
In general, geodesics are not optimal for all times. Instead, minimizers are geodesics by definition.

For 3D contact problem, we have that {\bf local minimizers are geodesics}. Indeed, first observe that the set of local minimizers is same if we consider the space $W^{1,1}$  or $W^{1,\infty}$, see Remark \ref{r-Linf}. Observe now that a local minimizer is a solution of the PMP (see \cite{agra-book}), and due to Proposition \ref{p-pmp} a curve is a solution of the PMP if and only if it is a geodesic. See more details in \cite{davidenote}.

A 3D contact manifold is said to be ``complete'' if all geodesics are defined for all times. This is the case for the problem \pmec\ defined in Section \ref{s-pmec} below. \b{In the following, for simplicity of notation, we always deal with complete 3D contact manifolds}.

In the following we denote by $(q(t),p(t))=e^{t\vec{H}}(q_{0},p_{0})$ the unique solution at time $t$ of  the Hamiltonian system
\bqn
\dot q=\partial_p H,~~~\dot p=-\partial_q H,
\eqnn
with initial condition $(q(0),p(0))=(q_{0},p_{0})$. Moreover we denote by $\pi:T^{*}M\to M$ the canonical projection $(q,p)\mapsto q$.
\bdeff
Let $(M,\textrm{span}\{X_{1},X_{2}\})$ be a 3D contact manifold and $q_0\in M$.
  Let
  $ \Lambda_{q_{0}}:=\{p_{0}\in T^{*}_{q_{0}}M | \, H(q_{0},p_{0})=1/2\}$.
We define the \emph{exponential map} starting from $q_{0}$ as
\bqn \label{eq:expmap}
\EXP_{q_{0}}: \Lambda_{q_{0}}\times \R^{+} \to M, \qquad \EXP_{q_{0}}(p_{0},t)= \pi(e^{t\vec{H}}(q_{0},p_{0})).
\eqn
\edeff

Next, we recall the definition of cut and conjugate time.
\bdeff \label{def:cut} Let $q_{0}\in M$ and
$\g$ be a geodesic parameterized by \srarc\  starting from $q_{0}$.
The \emph{cut time} for $\g$ is $\tcut(\g)=\sup\{t>0,\ |\ \g|_{[0,t]} \text{ is optimal}\}$. The \emph{cut locus} from $q_{0}$ is the set
$$\Cut(q_{0})=\{q(\tcut(\g))\ |\ \g \text{ geodesic parameterized by \srarc\  starting from }q_{0}\}.$$
\edeff

\bdeff \label{def:con} Let $q_{0}\in M$ and
$\g$  be a geodesic parameterized by \srarc\  starting from $q_{0}$ with initial covector $p_{0}$.
The \emph{first conjugate time} of $\g$ is
$$\tconj(\g)=\min\{t>0\ |\  (p_{0},t)  \text{ is a critical point of}\  \EXP_{q_{0}}\}.$$ The \emph{conjugate locus} from $q_{0}$ is the set $\Con(q_{0})=\{q(\tconj(\g))\  |\ \g \text{ \srarc\  geodesic from }q_{0}\}$.

\edeff
A geodesic loses its local optimality at its first conjugate locus. However a geodesic can lose optimality for ``global'' reasons. Hence we introduce the following:

\bdeff \label{def:max} Let $q_{0}\in M$ and
$\g$ be a geodesic parameterized by \srarc\  starting from $q_{0}$.
We say that $t_{max}>0$ is a \emph{Maxwell time} for $\g$ if there exists another   geodesic  $\bar q(.)$, parameterized by \srarc\  starting from $q_{0}$ such that $q(t_{max})=\bar q(t_{max})$
\edeff

It is well known that, for a geodesic $\g$,  the cut time $\tcut(\g)$ is either equal to the first conjugate time or to the first Maxwell time, 
see for instance \cite{agra-exp}.  Moreover, we have (see again \cite{agra-exp}):
\bt \label{corr:1}
Let $\gamma$ be a geodesic starting from $q_0$ and let $\tcut$ and $\tconj$ be its cut and conjugate times (possibly $+\infty$). Then
\bi
\iii $\tcut\leq\tconj$;
\iii $\gamma$ is globally optimal from $t=0$ to $\tcut$ and it is not globally optimal from $t=0$ to $\tcut+\eps$, for every $\eps>0$;
\iii $\gamma$ is locally optimal from $t=0$ to $\tconj$ and it is not locally optimal from $t=0$ to $\tconj+\eps$, for every $\eps>0$.\ei
\et

\brem
In 3D contact geometry (and more in general in sub-Riemannian geometry) the exponential map is never a local diffeomorphism in a  neighborhood of a point. As a consequence, spheres are never smooth and both the cut and the conjugate locus from $q_{0}$ are adjacent to the point $q_{0}$ itself, i.e. $q_0$ is contained in their closure (see \cite{agratorino}).
\erem

\section{Definition and study of \pmec}
\label{s-pmec}

In this section we introduce the auxiliary mechanical problem \pmec. The study of solutions of such problem is the main tool that we use to prove Theorem \ref{t-maini}.

We first define the mechanical problem \pmec.
\bd
\iii{\pmec:}
Fix $\xi>0$ and  $(x_{in},y_{in},\theta_{in}),(x_{fin},y_{fin},\theta_{fin})\in \R^2\times S^1$.
In the space of $L^1$ controls $u(\cdot),v(\cdot):[0,T]\to\R$,
find the solutions of:
\bqn
&&(x(0),y(0),\theta(0))=(x_{in},y_{in},\theta_{in}),~~~(x(T),y(T),\theta(T))=(x_{fin},y_{fin},\theta_{fin}) ,\nn\\
&&\Pt{\ba{c}
\frac{d x}{d\tau}(\tau)\\
\frac{d y}{d\tau}(\tau)\\
\frac{d \th}{d\tau}(\tau)
\ea
}=u(\tau)\Pt{\ba{c}
\cos(\th(\tau))\\
\sin(\th(\tau))\\
0
\ea}+v(\tau)\Pt{\ba{c}
0\\
0\\
1
\ea}\nn\\
&&   \int_0^T\sqrt{\xi^2u(\tau)^2+v(\tau)^2}~d\tau\to\min ~~~(\mbox{here $T$ is free})\label{eq-KOST}
\eqn
\ed

This problem (which cannot be interpreted as a problem of reconstruction of planar curves, as explained in \cite{nostro-vision}) has been completely solved in a series of papers by one of the authors  (see \cite{yuri1,yuri2,yuri3}).

\begin{remark} \label{remark:xi}
Observe that \pmec\ (as well as \pprojective\ and \pcurve) depend on a parameter $\xi>0$. It is easy to reduce our study to the case $\xi=1$. Indeed, consider the problem \pmec\ with a fixed $\xi>0$, that we call \pmec $(\xi)$. Given a curve $\g$ with cost $C_\xi$, apply the dilation $(x,y)\to(\frac{1}\xi x, \frac{1}\xi y)$ to find a curve $\tilde q(.)$. This curve has boundary conditions that are dilations of the previous boundary conditions, and it satisfies the dynamics for \pmec. If one considers now its cost $C_1$ for the problem \pmec $(1)$, one finds that $C_1= C_\xi$. Hence, the problem of minimization for all \pmec\ is equivalent to the case \pmec $(1)$. The same holds for \pprojective, \pcurve, with an identical proof. For this reason, we will fix $\xi=1$ from now on.
\end{remark}

Remark that \pmec\ is a 3D contact problem. Then, one can use the techniques given in Section \ref{s-sR} to compute the minimizers. This is the goal of the next section.

\subsection{Computation of geodesics for \pmec}
\label{s-geo-pmec}

In this section, we compute the geodesics for \pmec\ with $\xi=1$, and prove some properties that will be useful in the following. First observe that for \pmec\ there is existence of minimizers for every pair $(q_{in},q_{fin})$ of initial and final conditions, and minimizers are geodesics. See \cite{yuri1,yuri2,yuri3}. Moreover, geodesics are analytic, see Remark \ref{analytic}.

Since \pmec\ is 3D contact, we can apply Proposition \ref{p-pmp} to compute geodesics. We recall that we have
 $$
 M=\R^2\times S^1,~~
q=(x,y,\th),~~p=(p_1,p_{2},p_{3}),~~X_{1}=\left(\ba{c}
\cos(\th)\\
\sin(\th)\\0
\ea
\right),~~~
X_{2}=\left(\ba{c}
0\\
0\\1
\ea
\right).
$$
Hence, by Proposition \ref{p-pmp}, we have
$$H=\frac12\left(  (p_1\cos\th+p_2\sin\th)^2+p_3^2\right),$$
and the Hamiltonian equations are:
\bqn
&\dot x=\frac{\partial H}{\partial p_1}=h(q,p)\cos\th,&\dot p_1=-\frac{\partial H}{\partial x}=0,\nn\\
&\dot y=\frac{\partial H}{\partial p_2}=h(q,p)\sin\th,&\dot p_2=-\frac{\partial H}{\partial y}=0,\nn\\
&\dot \th=\frac{\partial H}{\partial p_3}=p_3,&\dot p_3=-\frac{\partial H}{\partial \th}=-h(q,p)(-p_1 \sin \th+ p_2 \cos\th),\nn
\eqn
where $h(q,p)=p_1\cos\th+p_2\sin\th$. Notice that this Hamiltonian system is integrable in the sense of Liouville, since we have enough constant of the motions in involution. Moreover, it can be solved easily in terms of elliptic functions. Setting $p_1=P_r\cos P_a$, $p_2= P_r\sin P_a$ one has $h(p,q)=P_r\cos(\th-P_a)$ and $\th(t)$ is solution of the pendulum like equation $\ddot\th=\frac12P_r^2\sin(2(\th-P_a))$.  Due to invariance by rototranslations, the initial condition on the $q$ variable can be fixed to be $(x_{in},y_{in},\th_{in})=(0,0,0)$, without loss of generality. The initial condition on the $p$ variable is such that $H(0)=1/2$. Hence $p(0)$ must belong to the cilinder
\begin{equation}\label{cilinder}
C=\{(p_1,p_2,p_3)~|~ p_1^2+p_3^2=1\}.
\end{equation}
In the following we use the notation of \cite{yuri1,yuri2,yuri3}. Introduce coordinates $(\nu, c)$ on $C$ as follows:
\begin{equation} \label{betac}
\sin (\nu/2) = p_1 \cos \th +  p_2\sin \th, \qquad
\cos (\nu/2) = - p_3, \qquad  c = 2 (p_2 \cos \th - p_1 \sin \th),
\end{equation}
with $\nu=2\th +\pi\in 2 S^1$. Here $2 S^1 = \R /(4 \pi \Z)$ is   the double covering  of the standard circle $S^1 = \R /(2 \pi \Z)$.\\

In these coordinates, the Hamiltonian system reads as follows:
\bqn
&\label{ham_vert}
\dot \nu = c, \quad \dot c = -\sin \nu, \qquad (\nu, c) \in   (2 S^1_{\nu}) \times \R_c, \\
&\label{ham_hor}
\dot x = \sin \frac{\nu}{2} \cos \theta, \quad \dot y = \sin \frac{\nu}{2} \sin \theta, \quad \dot \theta = - \cos \frac{\nu}{2}.
\eqn

Note that the curvature of the curve $(x(.), y(.))$ is equal to
\bqn
K = \frac{\dot x \ddot y - \ddot x \dot y}{(\dot x ^2 + \dot y^2)^{3/2}} = - \cot (\nu/2).\label{eq-k}
\eqn
We now define cusps for geodesics of \pmec. Recall that both the geodesics and the corresponding controls are analytic.
 \bdeff
 Let $\g=(x(\cdot),y(\cdot),\th(\cdot))$ be a geodesic of \pmec, parameterized by sR-arclength. We say that $\tcusp$ is a cusp time for $\g$ (and
 $q(\tcusp)$ a cusp point) if $u(.)$ changes its sign at $\tcusp$. We say that the restriction of  $q(\cdot)$ to an interval $[0,T]$ has no {\bf internal cusps} if no  $t\in]0,T[$ is a cusp time.
 \edeff

Given a curve $\g$ with a cusp point at $\tcusp$, we have that its projection on the plane $x(\cdot),y(\cdot)$ has a planar cusp at $\tcusp$ as well, see Figure \ref{fig:f-cusp}. More precisely, we have the following lemma.
\bl \label{l-cusps} A geodesic $\gamma$ (without angular cusps) has a cusp at $\tcusp$ if and only if $\lim_{t\to\tcusp} |K(t)|=\infty$.
\el
\bproof First observe that $\gamma$ has an internal cusp at $\tcusp$ if, for $t\to\tcusp$, it holds $u(t)\to 0$ and $v(t) \not\to 0$, i.e. using \r{ham_hor} one has $u(t)=\sin(\nu/2)\to 0 $ and $v(t)=-\cos(\nu/2)\not\to 0$. This is equivalent to $K(t)=- \cot (\nu/2)\to \infty$, by using \r{eq-k}.
\eproof

Also observe that one can recover inflection points of the planar curve $x(\cdot),y(\cdot)$ from the expression of $\g$. Indeed, at an inflection point of the planar curve, we have that the corresponding $\g$ satisfies $K=0$ and $\nu=\pi+2\pi n$, with $n\in\Z$.

\subsection{Qualitative form of the geodesics}
\label{s-qualitative}
Equation (\ref{ham_vert}) is the  pendulum equation
\bqn
\ddot\nu=-\sin\nu, \qquad \nu \in 2 S^1 = \R /(4 \pi \Z),
\label{eq-pend}
\eqn
whose phase portrait is shown in Figure~\ref{fig:phase}.

\immagine[8]{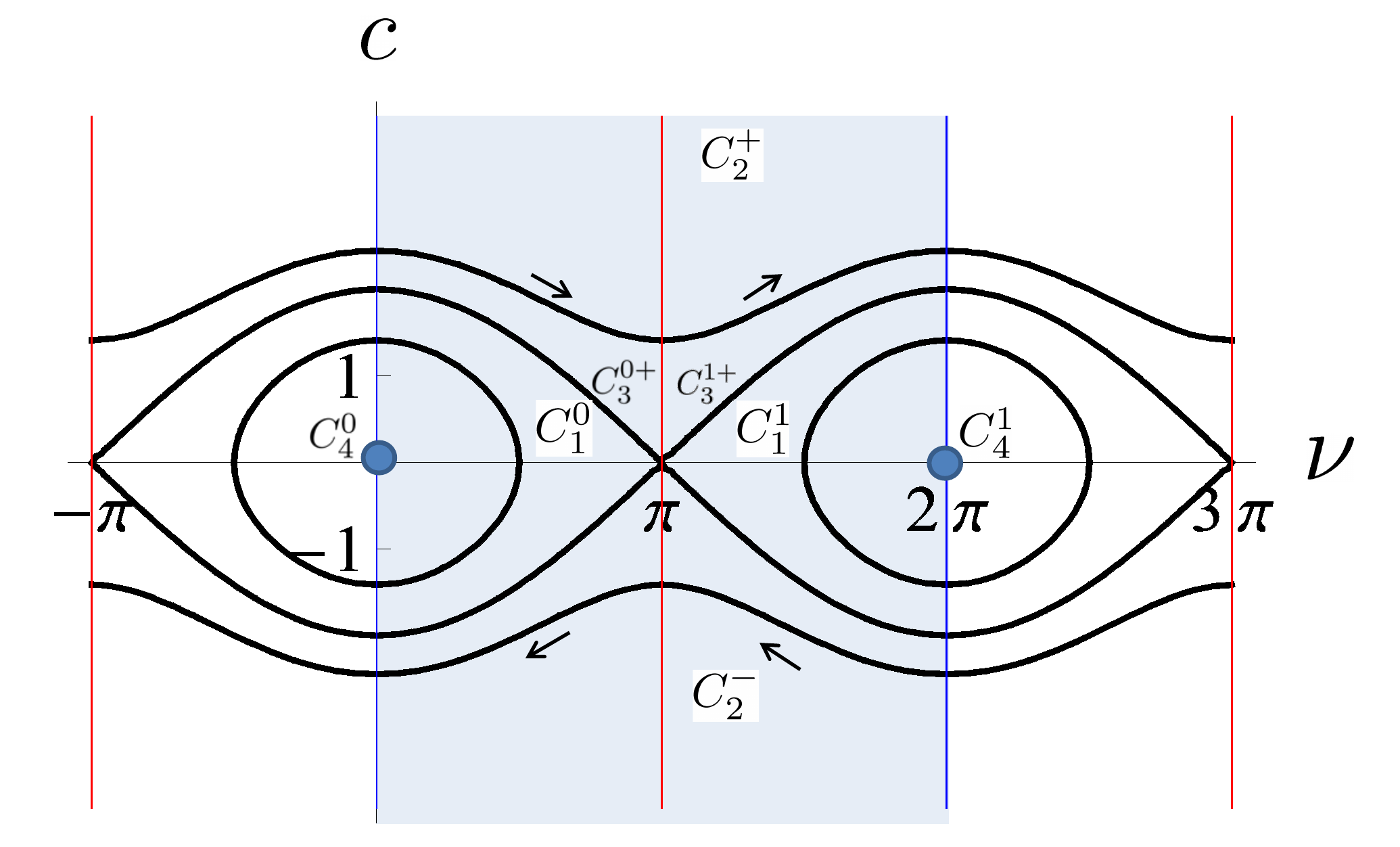}{Phase portraits of the pendulum equation, with the 5 types of trajectories.}

\begin{figure}[htbp]
\includegraphics[width=0.32\textwidth]{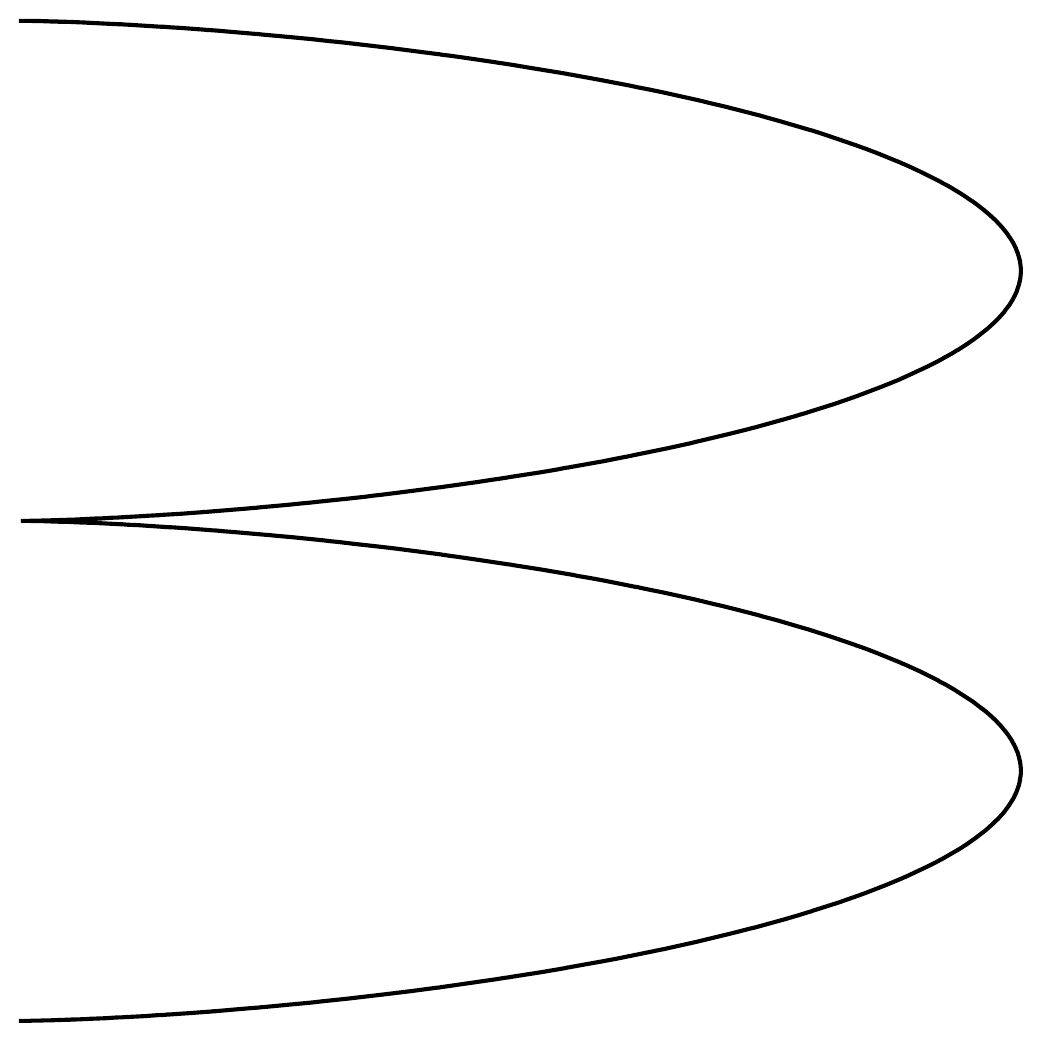}
\hfill
\includegraphics[width=0.32\textwidth]{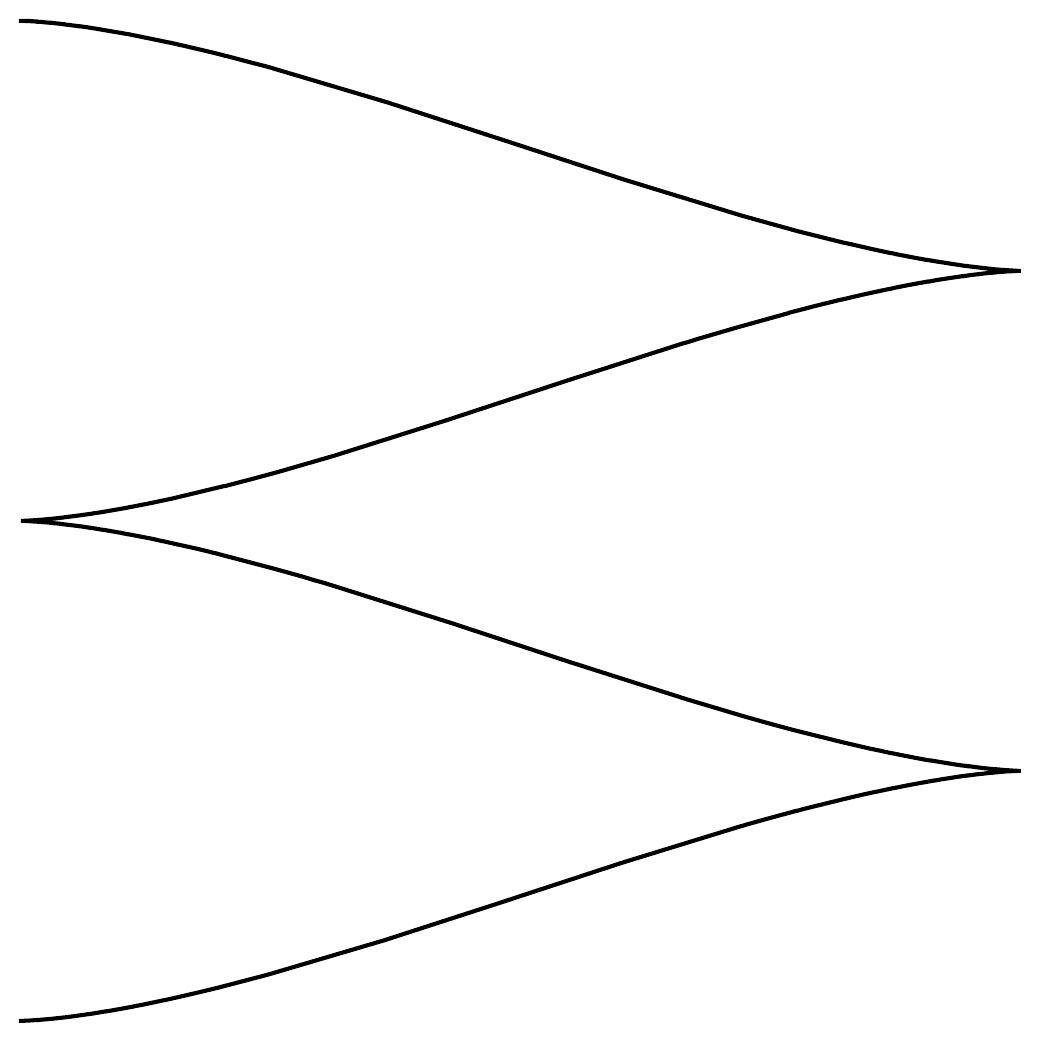}
\hfill
\includegraphics[width=0.32\textwidth]{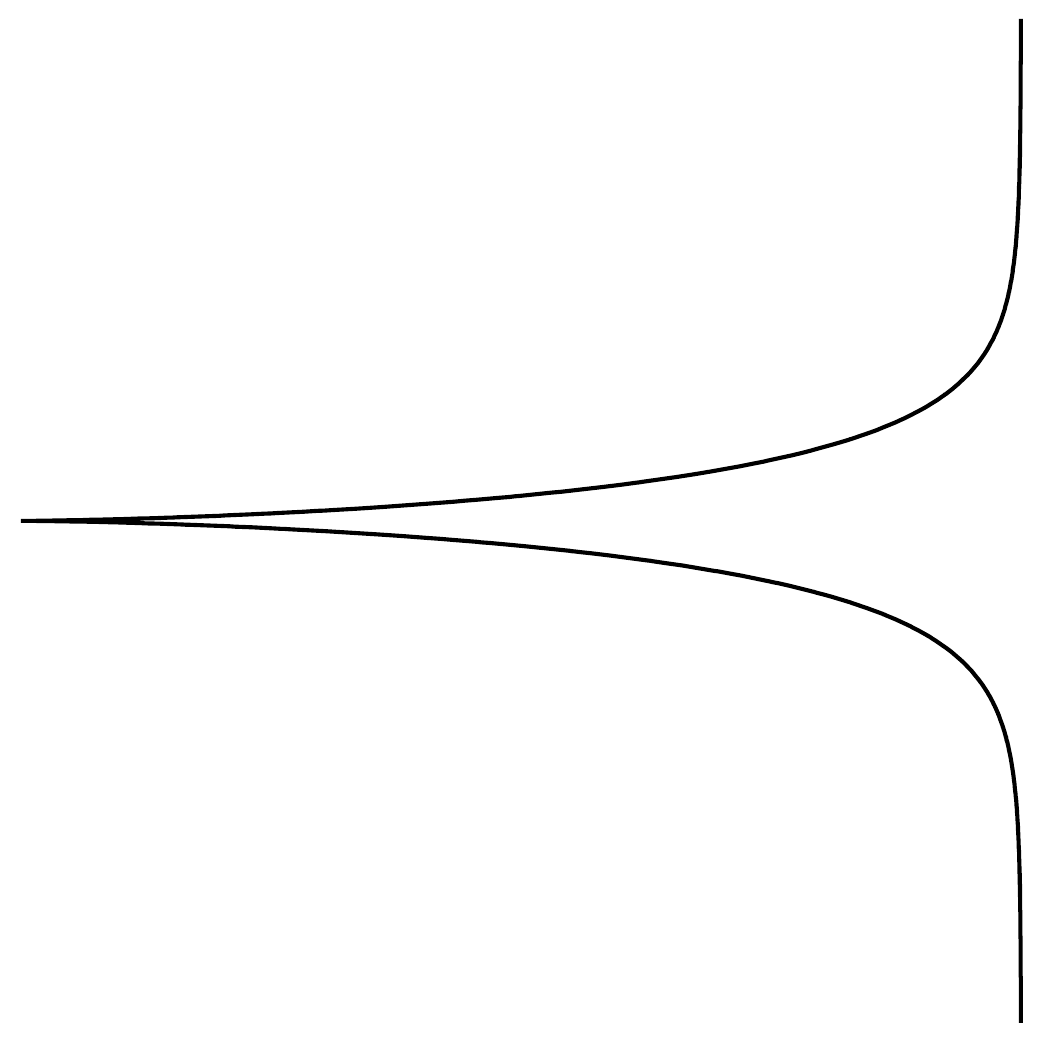}
\\
\parbox[t]{0.31\textwidth}
{\caption{Trajectory of type R.}\label{fig:xyC1}}
\hfill
\parbox[t]{0.32\textwidth}
{\caption{Trajectory of type O.}\label{fig:xyC2}}
\hfill
\parbox[t]{0.34\textwidth}
{\caption{Trajectory of type Sep.}\label{fig:xyC3}}
\end{figure}

 There exist 5 types of geodesics corresponding the different pendulum trajectories.
\begin{enumerate}
\item
Type S:  stable equilibrium of the pendulum: $\nu\equiv0$. For the corresponding planar trajectory, in this case we have
$(x(t), y(t)) \equiv (0, 0)$. These are the only geodesics with  angular cusps.
\item Type U:  unstable equilibria of the pendulum:  $\nu\equiv \pi$  or $\nu\equiv-\pi$.  For the corresponding planar trajectory, in this case we have
$(x(t), y(t)) = (t, 0)$ or $(x(t), y(t)) = (-t, 0)$, i.e. we get a  straight line.
\item Type R: rotating pendulum.  For the corresponding planar trajectory, in this case we have that
$(x(t), y(t))$ has infinite number of cusps and no inflection points (Fig.~\ref{fig:xyC1}). Note that in this case $\theta$ is a monotone function.
\item Type O: oscillating pendulum. For the corresponding planar trajectory, in this case we have that  $(x(t), y(t))$ has infinite number of cusps and infinite number of inflection points (Fig.~\ref{fig:xyC2}). Observe that between two cusps we have an inflection point, and between two inflection points we have a cusp.
\item Type Sep: separating trajectory of the pendulum. For the corresponding planar trajectory, in this case we have that
$(x(t), y(t))$ has one cusps and no inflection points (Fig.~\ref{fig:xyC3}).
\end{enumerate}
The explicit expression of geodesics in terms of elliptic functions are recalled in Appendix \ref{app:A}.

Recall that, for trajectories of type R, O and Sep, the cusp occurs whenever $\nu(t)=2\pi n$, with $n\in \Z$, since in this case one has from Lemma \ref{l-cusps} that $K(t)\to\infty$ for $t\to\tcusp$.

 \subsection{Optimality of geodesics}
 \renewcommand{\SS}{\mathbb{S}}
 \newcommand{\TT}{\mathbb{T}}
 
\begin{figure}[htbp]
\includegraphics[height=5.cm]{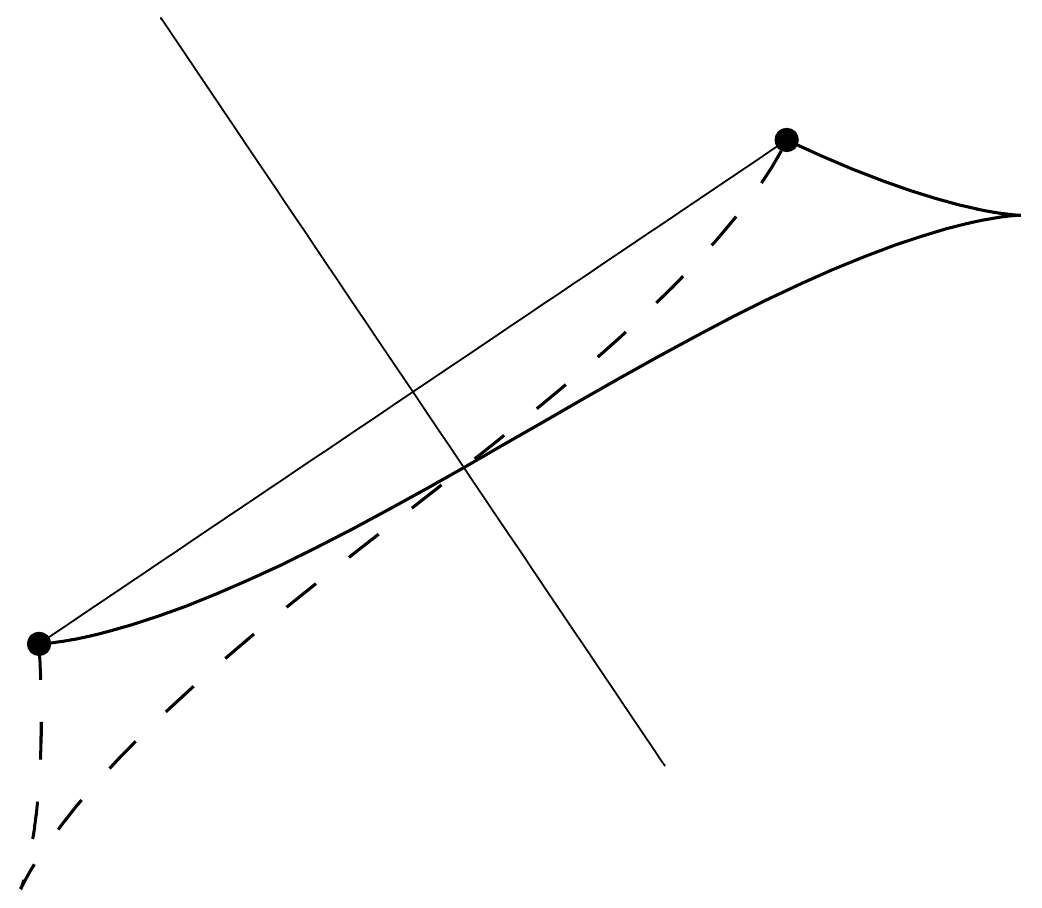}
\hfill
\includegraphics[height=3.cm]{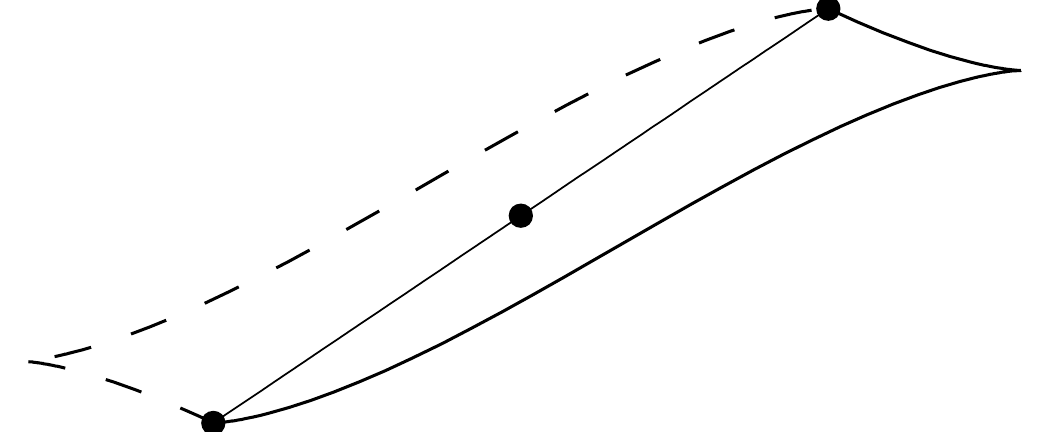}
\\
\parbox[t]{0.48\textwidth}{\caption{Action of   $\SS$ on $t \mapsto (x(t), y(t))$.}\label{fig:eps2}}
\hfill
\parbox[t]{0.48\textwidth}{\caption{Action of $\TT$   on $t \mapsto (x(t), y(t))$.}\label{fig:eps5}}
\end{figure}

Let $q(.)=(x(.),y(.),\th(.))$ be a geodesic parameterized by sub-Riemannian arclength $t \in [0, T]$. Consider the following two mappings of geodesics\footnote{Such mappings are denoted by $\eps^2,\eps^5$ in \cite{yuri1,yuri2,yuri3}, respectively.}:
$$
 \SS,\TT : \ q(.) \mapsto q_\SS(.),q_\TT(.),\qquad \mbox{~~~with~~}q(.):[0,T]\to \R^2\times S^1
$$
where
\bqn
\th_\SS(t) &=& \th(T) - \th({T-t}), \nn\\
 x_\SS(t) &=& -\cos \th(T) (x(T) - x({T-t})) - \sin \th(T) (y(T) - y({T-t})), \nn\\
 y_\SS(t) &=& -\sin \th(T) (x(T) - x({T-t})) + \cos \th(T) (y(T) - y({T-t})),
\eqnn
and
\bqn
 \th_\TT(t) &=& \th({T-t}) -  \th({T}), \nn\\
 x_\TT(t) &=& \cos \th(T) (x({T-t}) - x({T})) + \sin \th(T) (y({T-t}) - y({T})), \nn\\
 y_\TT(t) &=& -\sin \th(T) (x({T-t}) - x({T})) + \cos \th(T) (y({T-t}) - y({T})).
\eqnn
Modulo rotations of the plane $(x,y)$, the mapping $\SS$ acts as reflection of the curve $(x(.), y(.))$  in the middle perpendicular to the segment that connects the points $(x(0), y(0))$ and $(x(T), y(T))$; the mapping  $\TT$ acts as reflection in the midpoint of this segment. See Figures \ref{fig:eps2} and \ref{fig:eps5}.

\begin{figure}[htbp]
\includegraphics[height=5cm]{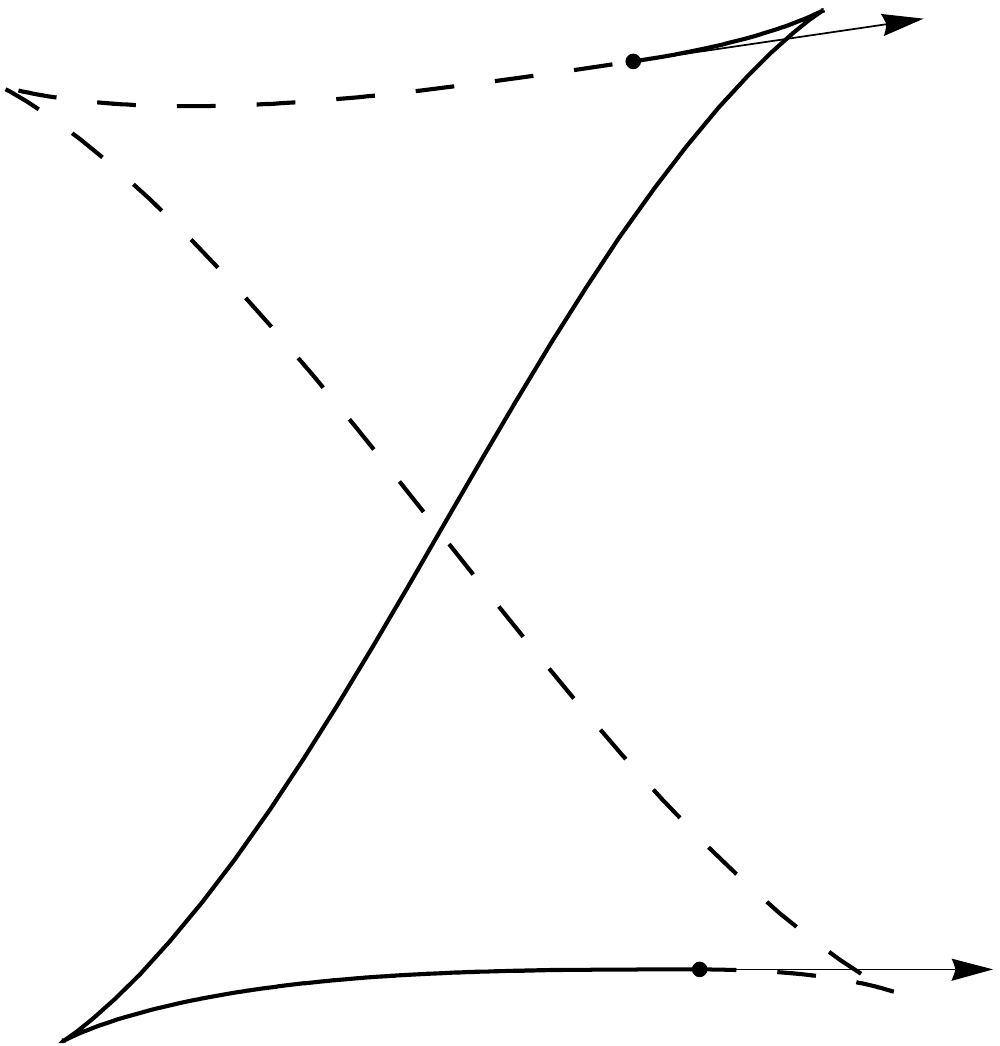}
\hfill
\includegraphics[height=5cm]{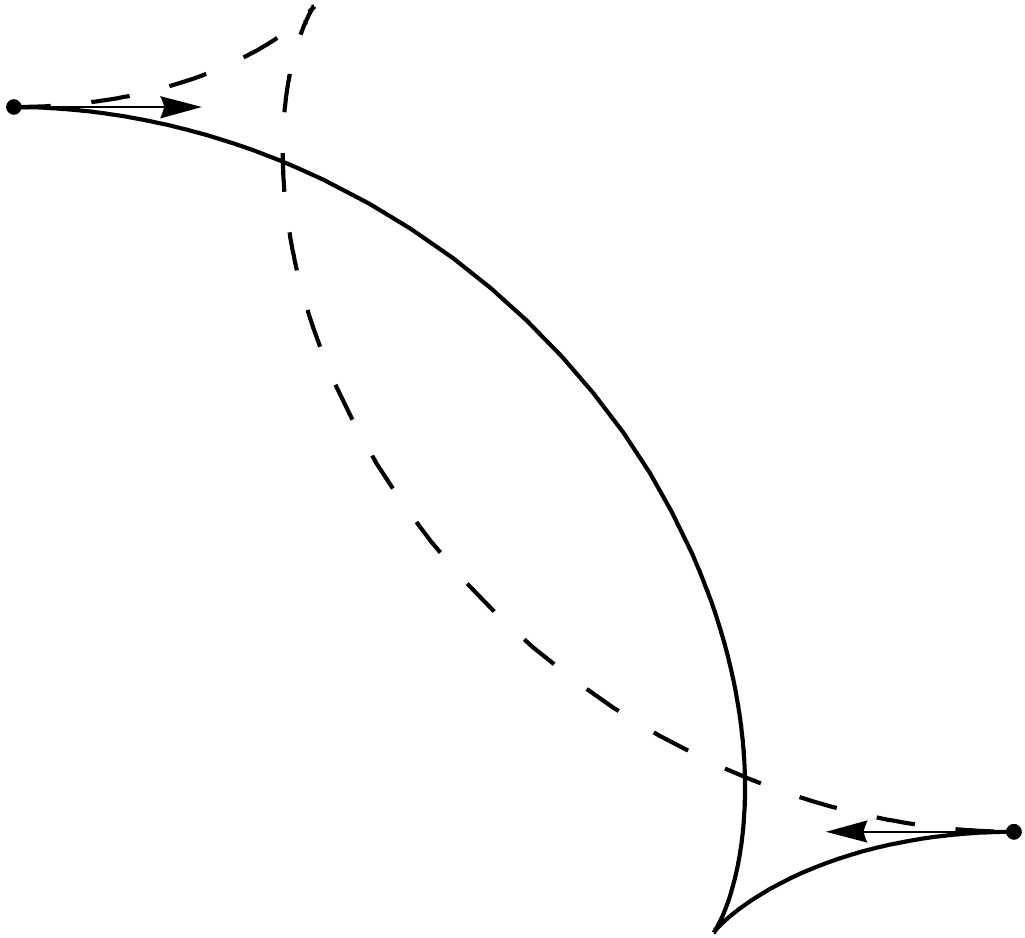}
\\
\parbox[t]{0.45\textwidth}{\caption{Maxwell point for reflection $\SS$.}\label{fig:Max2}}
\hfill
\parbox[t]{0.45\textwidth}{\caption{Maxwell point for reflection $\TT$.}\label{fig:Max5}}
\end{figure}

A point $q(t)$ of a trajectory $\g$ is called a Maxwell point corresponding to the reflection $\SS$ if $q(t) = q_\SS(t)$ and $q(\cdot) \not\equiv q_\SS(\cdot)$. The same definition can be given for $\TT$. Examples of Maxwell points for the reflections $\SS$ and $\TT$ are shown at Figures \ref{fig:Max2}  and~\ref{fig:Max5}.

The following theorem  proved in \cite{yuri1,yuri2,yuri3} describes optimality of geodesics.

\bt
A geodesic $\g$ on the interval $[0, T]$,  is optimal if and only if
each point $q(t)$, $t \in (0, T)$, is neither a
 Maxwell points corresponding to   $\SS$ or $\TT$, nor the limit of a sequence of Maxwell points.
\et
Notice that if a point $q(t)$ is a limit of Maxwell points then it is a Maxwell point or a conjugate point.

Denote by $\tpend$ the period of motion of the pendulum~(\ref{eq-pend}). It was proved in~\cite{yuri3} that the cut time satisfies the following:
\begin{itemize}
\item
$\tcut = \frac 12 \tpend$ for geodesics of type R,
\item
$\tcut \in (\frac 12 \tpend,\tpend)$ for geodesics of type O,
\item
$\tcut = + \infty = \tpend$ for geodesics of types S, U and Sep.
\end{itemize}

\bc \label{corr:cuspcut}
Let $\g$ be a geodesic. Let $\tcusp,$ and $\tcut$ be the first cusp time and the cut time (possibly $+\infty$).Then $\tcusp\leq\tcut$.
\label{c-1}
\ec
\bproof For geodesics of types R and O it follows from the phase portrait of pendulum~(\ref{eq-pend}) that there exists $t \in (0,  \frac 12 \tpend)$ such that $\nu(t) = 2 \pi n$. This implies that $K(t)\to+\infty$, and, by Lemma \ref{l-cusps}, we have a cusp point for such $t$. Then $\tcusp \leq \frac  12 \tpend \leq \tcut$.

For geodesics of types S, U and Sep, the inequality  $\tcusp \leq \tcut$ is obvious since $\tcut = + \infty$.
\eproof

\bc
\label{nonapofantic}
Let $\g$ defined on $[0,T]$ be a minimizer having an internal cusp. Then any other minimizer between $q(0)$ and $q(T)$
has an internal cusp.
\ec
\bproof It was proved in~\cite{yuri3} that for any points $q_0, q_1   \in \R^2 \times S^1$,
there exist either one or two minimizers connecting $q_0$ to $q_1$. Moreover, if there are two  such minimizers $\g$ and $\tilde{q}(\cdot)$, then $\tilde{q}(\cdot)$ is obtained from $\g$ by a reflection $\SS$ or $\TT$. So if $\g$ has an internal cusp, then $\tilde{q}(\cdot)$ has an internal cusp as well.
\eproof

\section{Equivalence of problems}
\label{s-problems}

In this section, we state precisely the connections between minimizers of problems  \pcurve, \pprojective\ and \pmec\ defined above. 
The problems are recalled in Table 1 for the reader's convenience. We also prove that there exists minimizers of \pcurve\ that are absolutely continuous but not Lispchitz.

\begin{center}
\begin{tabular}{|l|}
\hline
\\
{\bf Notation} \\\\
$\displaystyle
q=\left(\ba{c}
x\\
 y\\
\th
\ea\right),~~
X_{1}=\left(\ba{c}
\cos\th\\
 \sin\th\\
0
\ea\right),
~~
X_{2}=\left(\ba{c}
0\\
0\\
1
\ea\right),
$\\\\
here $\ul{x}:=(x,y)\in\R^2$ and $\theta\in S^1$ or $P^1$ as specified below. We denote with $s$ the plane-arclength\\
 parameter and with $t$ the \srarc\ parameter. In all problems written below we have the following:\\\\
$\bullet$ initial and final conditions  $(x_{in},y_{in},\theta_{in}),(x_{fin},y_{fin},\theta_{fin})$ are given;\\
$\bullet$ the final time $T$ (or length $\ell$) is free.\\\hline

{\bf Problem} \pcurve:\\
$\displaystyle
q\in \R^2\times S^1
~~~\dot q=X_{1}+v X_{2},~~\int_0^\ell \sqrt{\xi^2+ v^2}~ds=\int_0^\ell \sqrt{\xi^2+K(s)^2} ds\to\min
$\\\\\hline

{\bf Problem} \pmec:\\
$\displaystyle
q\in \R^2\times S^1
~~~~\dot q=uX_{1}+v X_{2},~~\int_0^T \sqrt{\xi^2 u^2+ v^2}~d t\to\min
$\\\\\hline
{\bf Problem} \pprojective:\\
$\displaystyle
q\in \R^2\times P^1
~~~~\dot q=uX_{1}+v X_{2},~~\int_0^T \sqrt{\xi^2 u^2+ v^2}\,dt=\int_0^T \| \dot {\bf x} \| \sqrt{\xi^2+ K(t)^2}\,dt\to\min
$
\\\\\hline
\end{tabular}
 Table 1. The different problems we study in the paper.
\end{center}

First notice that the problem \pmec\ admits a solution, as shown in \cite{yuri1,yuri2,yuri3}. The same arguments apply to \pprojective, for which existence of a solution is verified as well, see \cite{suzdal}.

Also recall that the definitions of \pprojective\ and \pmec\ are very similar, with the only difference that $\theta\in P^1$ or $\theta\in S^1$, respectively. This is based on the fact that $\R^2\times S^1$ is a double covering of $\R^2\times P^1$. Moreover, both the dynamics and the infinitesimal cost in \pmec\ are compatible with the projection $\R^2\times S^1 \to \R^2\times P^1$. Thus, the geodesics for \pprojective\ are the projection of the geodesics for \pmec. Then, locally the two problems are equivalent. If we look for the minimizer  for  \pprojective\ from $(x_{in},y_{in},\th_{in})$ to $(x_{fin},y_{fin},\th_{fin})$, then it is the shortest minimizer between the minimizing geodesics for  \pmec:
\bd
\i[minimizing geodesic $q_1(.)$]: connecting $(x_{in},y_{in},\th_{in})$ to $(x_{fin},y_{fin},\th_{fin})$;
\i[minimizing geodesic $q_2(.)$]: connecting $(x_{in},y_{in},\th_{in}+\pi)$ to $(x_{fin},y_{fin},\th_{fin})$;
\i[minimizing geodesic $q_3(.)$]: connecting $(x_{in},y_{in},\th_{in})$ to $(x_{fin},y_{fin},\th_{fin}+\pi)$;
\i[minimizing geodesic $q_4(.)$]: connecting $(x_{in},y_{in},\th_{in}+\pi)$ to $(x_{fin},y_{fin},\th_{fin}+\pi)$;
\ed
In reality, these four minizing geodesics are coupled two by two: indeed, $q_1$ and $q_4$ are geometrically the same curve, as well as $q_2$ and $q_3$. This is a direct consequence of the fact that one can reparametrize a geodesic backward in time, and as a consequence boudary conditions are transformed from $\th$ to $\th+\pi$. More precisely, there exists the following symmetry of geodesics for \pmec:
$(x(t), y(t), \theta(t)) \mapsto (x(t), y(t), \theta(t) + \pi),$
by replacing  $(u(t), v(t))$ with $(-u(t), v(t))$. See Figure \ref{fig:2sol}.

\begin{figure}[htb]
\begin{center}
\begin{lpic}{2sol}
\lbl{23,9;$(x_{in},y_{in})$}
\lbl{31,19,23;$\th_{in}$}
\lbl{10,4,23;$\th_{in}+\pi$}
\lbl{81,10;$(x_{fin},y_{fin})$}
\lbl{76,18,-31;$\th_{fin}$}
\lbl{92,1,-31;$\th_{fin}+\pi$}
\lbl{45,28,37;$q_1(.)\simeq q_4(.)$}
\lbl{54,15;$q_2(.)\simeq q_3(.)$}
\end{lpic}
\caption{Minimizing geodesic for \pprojective\ from minimizing geodesics for \pmec.}
\label{fig:2sol}
\end{center}
\end{figure}


It is also easy to prove that a minimizer of \pmec\ without cusps is also a minimizer of \pcurve. Indeed, take a minimizer $\g$ of \pmec\ without cusps, thus with $\|\dot \x(\tau)\|>0$ for $\tau\in\Pq{0,T}$. Then, reparametrize the time to have a spatial arclength parametrization, i.e. $u=\|\dot \x\| \equiv 1$ (this is possible exactly because it has no cusps). This new parametrization of $\g$ satisfies the dynamics for \pcurve\ and the boundary conditions. Assume now by contradiction that there exists a curve $\tilde q(.)$ satisfying the dynamics  for \pcurve\ and the boundary conditions with a cost that is smaller that the cost for $\g$. Then $\tilde q(.)$ also satisfies the dynamics for \pmec\ and boundary conditions, with a smaller cost, hence $\g$ is not a minimizer. Contradiction.

\subsection{Connection between curves of \pcurve\ and \pmec}

\label{s-genpmp}

In this section, we study in more detail the connection between curves of \pcurve\ and \pmec. First of all, observe that \pcurve\ and \pmec\ are defined on the same manifold $SE(2)$. Moreover, each curve $\Gamma(.)=(x(.),y(.),\th(.))$ satisfying the dynamics for \pcurve\ with a certain control $v(.)$, also satisfies the dynammics for \pmec\ with controls $u(.)\equiv 1$ and $v(.)$. For simplicity of notation, we give the following definition.

\bdeff Let $\Gamma(.)=(x(.),y(.),\th(.))$ be a curve in  $SE(2)$ satisfying  the dynamics for \pcurve\ with a certain control $v(.)$. We define the \b{corresponding curve} $q(.)$ for \pmec\ as the same parametrized curve $(x(.),y(.),\th(.))$, and the \b{corresponding pair} as the pair trajectory-control $(q(.),(u(.),v(.)))$ with $u(.)\equiv 1$.

We define the \b{corresponding reparametrized pair} $(q_1(.),(u_1(.),v_1(.)))$ for \pmec\ the time-reparametrization of the corresponding pair $(q(.),(u(.),v(.)))$ by \srarc, and the \b{corresponding reparametrized curve} as the curve $q_1(.)$.
\edeff

Recall that the time-reparametrization by \srarc\ of an admissible curve for \pmec\ is always possible. A detailed explanation for time-reparametrization of a curve with controls in $L^1$ to have controls in $L^\infty$ is given in \cite[Section 2.1.1]{suzdal}.

We now focus on solutions of the Pontryagin Maximum Principle (PMP). For \pcurve, one cannot apply the standard PMP since one cannot guarantee a priori that optimal controls are in $L^\infty$. For this reason, we apply a generalized version of the PMP which holds for $L^1$ controls (see  \cite[Thm 8.2.1]{vinter}). We have the following result.
\bt \label{t-equiv} Let $\Gamma(.)$ be a solution of the generalized PMP for \pcurve. Then the corresponding reparametrized 
curve is a solution of the standard PMP for \pmec.
\et

The proof of this Theorem is given in Appendix \ref{a-pmp}. Here we recall the
main steps of the proof:
\bd
\i[STEP 1:] we prove that if $(\Gamma(.), v(.))$ is a solution of the
generalized Pontryagin Maximum Principle for
\pcurve, then, the corresponding pair $(q(.),(u(.),v(.)))$ for \pmec\ is a solution of the generalized Pontryagin Maximum
Principle.
\i[STEP 2:] we prove that the corresponding reparametrized pair is a solution of the standard PMP.
\ed

We are now ready to discuss the connection between geodesics for \pcurve\ and \pmec.
\bp \label{p-geodesic} Let $\Gamma(.)$ be a geodesic for \pcurve. Then the corresponding reparametrized curve is a geodesic for \pmec.
\ep
\bproof Let $\Gamma(.)$ be a geodesic for \pcurve. By definition, for every sufficiently small interval its restriction is a global minimizer. Then it is a solution of the generalized PMP for \pcurve. Hence, applying the previous Theorem \ref{t-equiv}, we have that the corresponding reparametrized trajectory is a solution of the standard PMP for \pmec. This implies that it is a geodesic for \pmec, due to Proposition \ref{p-pmp}.
\eproof

\subsection{\pcurve\ admits minimizers which are absolutely continuous but not Lipschitz}
\label{ch:lav}

We now show that the problem \pcurve\ exhibits an interesting phenomenon: there exist absolutely continuous minimizers that are not Lipschitz. Other examples are given in \cite{sarychev-torres}.

Consider a geodesic of \pmec\ defined on $[0,T]$ having no internal cusp and corresponding to controls $u(\cdot)$ and $v(\cdot)$.  From Corollary \ref{c-1} it follows that it is optimal. Assume now that this geodesic has a cusp at $T$. Then, by Lemma \ref{l-cusps}, we have that for $t\to T$ it holds $u(t)\to0$ and $\kappa(t)\to\infty$.  Notice that $\sqrt{u(\tau)^2+v(\tau)^2}$ is integrable on $[0,T]$, since its integral is exactly the Carnot-Caratheodory distance \r{e-distanza}, that is finite, see e.g. \cite{yuri1}. Since the cost of \pmec\ and \pcurve\ coincide, we have that $\int_0^\ell \sqrt{1+K^2(s)}\,ds$ is finite. In particular, $K(.)$ is a $L^1$ function that is not $L^\infty$. Reparametrize time to have an admissible curve $\Gamma(.)$ for \pcurve, with control $\tilde v(.)$. Since $\tilde v(s)$ coincides with $\kappa(s)$, then $\tilde v(.)$ is a $L^1$ function that is not $L^\infty$. This means moreover that the trajectory $\Gamma(.)$ for \pcurve\ has unbounded control and it is not Lipschitz.

This phenomenon is extremely interesting in optimal control. Indeed, direct application of standard techniques for the computation of local minimizers, such as the Pontryagin Maximum Principle, would provide local minimizers in the ``too small'' set of controls $L^\infty([0,T],\R)$. In other words, the absolutely continuous minimizers that are not Lipschitz are not detected by the Pontryagin Maximum Principle. For this reason, we were obliged to use the generalized PMP for \pcurve\ in Theorem \ref{t-equiv}.

Instead, the auxiliary problem \pmec\ does not present this phenomenon, since by re-parametrization one can always reduce to the set $L^\infty([0,T],\R)$.
\section{Existence of minimizing curves}
\label{s-main}

In this section we prove the main results of this paper, proving Theorem \ref{t-maini}. We characterize the set of boundary conditions for which a solution of \pcurve\ exists. We show that the set of boundary conditions for which a {\bf solution} exists coincides with the set of boundary conditions for which a {\bf local minimizer} exists. Moreover, such set coincides with the set of boundary conditions for which a {\bf geodesics} joining them exists. We also give some properties of such set.

After this theoretical result, we show explicitly the set of initial and final points for which a solution exists, computed numerically. For more details on this subject, see \cite{future}.

From the following result, Theorem \ref{t-maini} follows.

\bt[main result] \label{th:main}
  Fix an initial and a final condition $q_{in}=(x_{in},y_{in},\th_{in})$ and $q_{fin}=(x_{fin},y_{fin},\th_{fin})$ in $\R^2\times S^1$ . Let $\g$ be a  minimizer for the problem \pmec\ from $q_{in}$ to $q_{fin}$. The only two possible cases are:
\begin{enumerate}
\iii  $\g$ has neither internal cusps nor  angular cusps. Then  $\g$ is a solution for \pcurve\ from $q_{in}$ to $q_{fin}$.
\iii $\g$ has at least an internal cusp or an angular cusp. Then  \pcurve\ from $q_{in}$ to $q_{fin}$ does not admit neither a global nor a  local minimum nor a geodesic.
\end{enumerate}
\label{t-strong}
\et
\newcommand{\q}{\g}
\bproof We use the notation $\q$ to denote trajectories for \pmec, and $\Gamma(.)$ for trajectories for \pcurve. Recall the results of Section \ref{s-problems}. Given a $\Gamma(.)=(x(.),y(.),\theta(.))$ trajectory of \pcurve, this gives naturally a $\q=(x(.),y(.),\theta(.))$ trajectory of \pmec. On the converse, a $\q=(x(.),y(.),\theta(.))$ trajectory of \pmec\ without cusps gives naturally a $\Gamma(.)$ trajectory of \pcurve, after reparametrization.

Fix an initial and a final condition
$q_{in}=(x_{in}, y_{in}, \theta_{in})$ and $q_{fin} = (x_{fin}, y_{fin}, \theta_{fin})$.
Take a solution $\q$ of \pmec. If $\q$ has no cusps, then one can reparametrize time to have a curve $\Gamma(.)$ solution of \pcurve. If $\q$ has cusps at boundaries, then the same re-parametrization (that can be applied, as explained in Section~\ref{s-genpmp}) gives the corresponding $\Gamma(.)$, that is a solution of \pcurve. The first part is now proved.\\

We prove the second part by contradiction. If $\q$ has an internal cusp, then any other solution of \pmec\ from $q_{in}$ to $q_{fin}$ has an
internal cusp, as proved in Corollary~\ref{nonapofantic}. By contradiction, assume that there exists $\bar\Gamma(.)$, either a solution (i.e. a global minimizer) of \pcurve\ from $q_{in}$ to $q_{fin}$, or a local minimizer, or a geodesic. In the three cases, the corresponding reparametrized curve on $SE(2)$ of \pmec, that we denote by $\bar q_1(.)$, has no cusps.

We first study the case of geodesics. Let $\bar \Gamma(.)$ be a \textbf{geodesic} of \pcurve. Then $\bar q_1(.)$ is a geodesic of \pmec\ between the same boundary conditions of $\bar \Gamma(.)$, due to Proposition \ref{p-geodesic}. Then, two cases are possible:
\bi 
\i Let $\bar q_1(.)$ be a solution, i.e. a global minimizer, for \pmec. Then both $\q$ and $\bar q(.)$ are minimizers, one with cusps and the other without cusps. This yelds a contradiction with Corollary \ref{nonapofantic}.
\i Let $\bar q_1(.)$ be a geodesic for \pmec\ that is not a global minimizer. We denote with $\Pq{0,T}$ the time-interval of definition of $\bar q_1(.)$. Then there exists a cut time $t_{cut}<T$ for $\bar q_1(.)$. Then there exists a cusp time $t_{cusp}\leq t_{cut}<T$ for $\bar q_1(.)$, see Corollary \ref{corr:cuspcut}. Then $\bar q_1(.)$ has a cusp. Contradiction.
\ei
We have a contradiction in both cases. Thus, if $\q$ has a internal cusp, there exists no geodesic of \pcurve\ from $q_{in}$ to $q_{fin}$.

We now study the case of \b{local minimizers}. Let $\bar \Gamma(.)$ be a local minimizer for \pcurve. Then, it is a solution of the generalized Pontryagin Maximum Principle \cite[Thm 8.2.1]{vinter}. Applying Theorem \ref{t-equiv}, we have that the corresponding reparametrized curve $\bar q_1(.)$ is a solution of the standard Pontryagin Maximum Principle for \pmec, and then it is a geodesic by Proposition \ref{p-pmp}. Since $\bar \Gamma(.)$ has no cusps, then $\bar q_1(.)$ has no cusps either, thus it is a global minimizer. Then both $\q$ and $\bar q_1(.)$ are global minimizers, one with cusps and the other without cusps. This yelds a contradiction with Corollary \ref{nonapofantic}.

Since \b{global minimizers} are special cases of local minimizers, we have the result for global minimizers too.\\

If instead $\q$ has an angular cusp, then $(x_{in},y_{in})=(x_{fin},y_{fin})$, see Remark \ref{r-cuspdistribuita}. In this case, assume that there exists $\bar\Gamma(.)$ either a solution of \pcurve\ (i.e. a global minimizer), or a local minimizer, or a geodesic. In the three cases, the corresponding reparametrized trajectory of \pmec\ $\bar q_1(.)$ must be of Type S, since there are no other geodesics steering $q_{in}$ to $q_{fin}$ with $(x_{in},y_{in})=(x_{fin},y_{fin})$. By construction, the solution of \pcurve\ is $\bar\Gamma(.)=(x_{in},y_{in},\th(.))$. Observing the dynamics for \pcurve\ in \r{e-pcurve}, one has that $x,y$ constant implies that the planar length is $\ell=0$, then we must have $\th_{in}=\th_{fin}$.
\eproof

\brem Observe that, as a corollary, we have proved that global minimizers, local minimizers and geodesics for \pcurve\ coincide.
\erem

\brem The last part of the proof has its practical interest. It shows the non-existence of a solution of \pcurve\ in the case of $(x_{in},y_{in})=(x_{fin},y_{fin})$. This means that, under this condition, it is possible to construct a sequence of planar curves $\gamma^n(.)$, each steering $(x_{in},y_{in},\th_{in})$ to $(x_{in},y_{in},\th_{fin})$ and such that the sequence of the costs of $\gamma^n(.)$ converges to the infimum of the cost, but that the limit trajectory $\gamma^*(.)$ is a curve reduced to a point, for which the curvature $\kappa$ is not well-defined. See Figure \ref{fig:unpunto}.
\erem
\immagine[8]{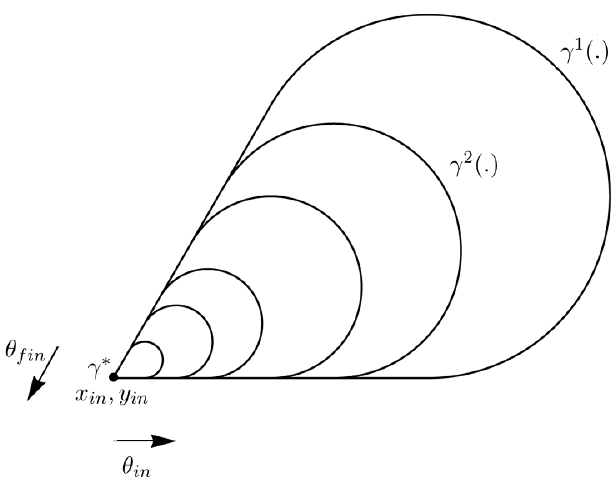}{Non-existence of a solution of \pcurve\ for $(x_{in},y_{in})=(x_{fin},y_{fin})$.}

\subsection{Characterization of the existence set}

In this section, we characterize the set of boundary conditions for which a solution of \pcurve\ exists, answering the second part of question {\bf Q2}. We recall that we just proved that the existence set does not change if we consider global or local minimizers or geodesics.

We prove here some simple topological properties of such set, and give some related numerical results.

\newcommand{\Sol}{\mathcal{S}}
\newcommand{\Lam}{\Lambda}

\bp Let $\Sol\subseteq \R^2\times S^1$ be the set of final conditions $q_{fin}=(x_{fin},y_{fin},\th_{fin})$ for which a solution of \pcurve\ exists, starting from $e:=(0,0,0)$. We have that $\Sol$ is arc-connected and non-compact.
\ep
\bproof For arc-connectedness, let $q_a,q_b\in \Sol$. This means that there exist two curves $q_1(.),q_2(.)$ steering $e$ to $q_a,q_b$, respectively. Then the concatenation of curves (with reversed time for $q_1(.)$) steers $q_a$ to $e$ to $q_b$. For non-compactness, observe that all points on the half-line $(t,0,0)$ are in $\Sol$.\eproof

Other properties of $\Sol$ (which are evident numerically\footnote{Formal proofs are given in \cite{future}.}) are the following:
\be
\i all points of $\Sol$ satisfy $x_{fin}\geq 0$;
\i if $q_{fin}\in\Sol$ satisfies $\theta_{fin}=\pi$, then it also satisfies $x_{fin}=0$; similarly, if $q_{fin}\in\Sol$ satisfies $x_{fin}=0$, then it also satisfies $\theta_{fin}=\pi$. The solutions of a problem with $q_{fin}=(0,y_{fin},\pi)$ have a cusp in $q_{fin}$.
\ee

\immagine[1.3]{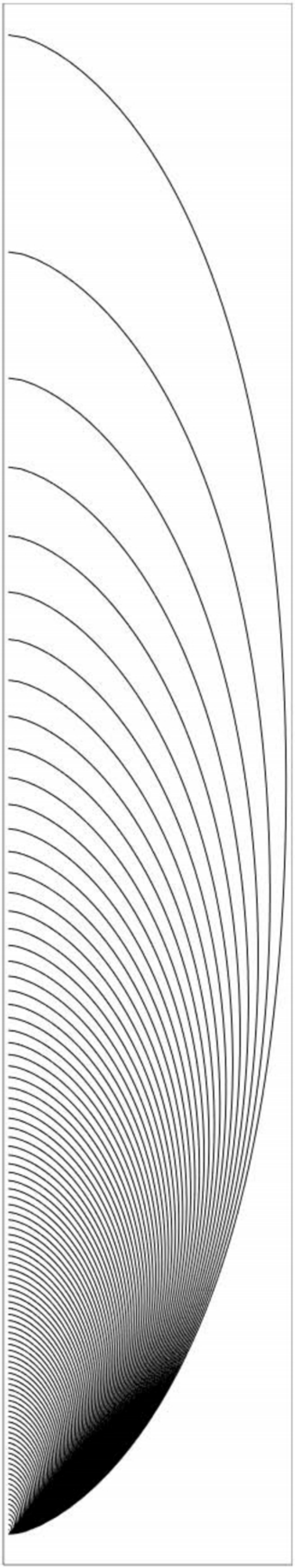}{Geodesics reaching $x=0$, upper plane.}


\brem The characterization of $\Sol$ is, in some sense, the continuation of the main results of the authors in \cite{suzdal}. There, we proved that there exist boundary conditions such that \pcurve\ did not admit a minimizer, i.e. that $\Sol$ is not the whole space $SE(2)$. Here we have described in bigger detail the set of boundary conditions $\Sol$ such that \pcurve\ admits a minimizer, together with proving that, given boundary conditions, the existence of a minimizer is equivalent to the existence of a local minimizer or a geodesic.
\erem

\section*{Acknowledgements}

The authors wish to thank Arpan~Ghosh and Tom~Dela~Haije, Eindhoven University of Technology, for the contribution with numerical computations and figures.

This research has been supported  by the European Research Council, ERC StG 2009 ``GeCoMethods'', contract number 239748, by the ANR ``GCM'', program ``Blanc--CSD'' project number NT09-504490, by the DIGITEO project ``CONGEO'', by Russian Foundation for Basic Research,  Project No.~12-01-00913-a, and by the Ministry of Education and Science of Russia within the federal program ``Scientific and Scientific-Pedagogical Personnel of Innovative Russia'', contract no. 8209.

 \appendix

\section{Explicit expression of geodesics in terms of elliptic functions}
\label{app:A}

In this section, we recall the explicit expressions of the geodesics for \pmec. They were first computed in \cite{yuri1}. 

The geodesics  are  expressed in sub-Riemannian arc-length $t$, and they are written in terms of Jacobian functions $\cn$, $\sn$, $\dn$, $\E$. For more details, see e.g. \cite{elliptic}. Here $(\nu,c)$ are the variables for the pendulum equation \r{ham_vert} and $(\varphi, k)$ are the corresponding action-angle coordinates that rectify its flow: $\dot \varphi = 1$,  $\dot k = 0$. See detailed explanations in \cite[Sec. 4]{yuri1}.

Since \pmec\ is invariant via rototranslations, we give geodesics starting from $(0,0,0)$ only.

Recall that we have classified geodesics of \pmec\ via the classification of trajectories of  the pendulum Eq.~(\ref{eq-pend}), see Section \ref{s-qualitative}. We have the following 5 cases.
\bi
\i The geodesic of type S has the simple expression $q(t)=(0,0,t)$. The projection on the plane gives the line reduced to the point $(0,0)$.
\i The geodesic of type U has the simple expression $q(t)=(t,0,0)$. The projection on the plane is the straight half-line $(t,0)$.
\i Geodesics of type R have the following expression :
\begin{align*}
&\cos \th(t) = \cn \varphi \cn (\varphi +t) + \sn \varphi \sn  (\varphi +t), \\
&\sin \th(t) = \sgn(\cos(\nu/2))(\sn \varphi \cn (\varphi +t) - \cn \varphi \sn (\varphi +t)), \\
&x(t) = \frac{\sgn(\cos(\nu/2))}{k} [ \cn \varphi (\dn \varphi - \dn (\varphi +t)) + \sn \varphi (t + \E(\varphi) - \E(\varphi +t))], \\
&y(t) = (1/k) [ \sn \varphi (\dn \varphi - \dn (\varphi +t)) - \cn \varphi (t + \E(\varphi) - \E(\varphi +t))].
\end{align*}

\i Geodesics of type O have the following expression :
\begin{align*}
&\cos \th(t) = k^2 \sn (\varphi /k) \sn (\varphi + t)/k + \dn (\varphi /k) \dn (\varphi + t)/k, \\
&\sin \th(t) = k(\sn (\varphi /k) \dn (\varphi + t)/k - \dn (\varphi /k) \sn (\varphi + t)/k), \\
&x(t) = \sgn(c) k [\dn (\varphi /k)(\cn (\varphi /k) - \cn (\varphi + t)/k) + \sn (\varphi /k) (t/k + \E(\varphi /k) - \E((\varphi + t)/k)], \\
&y(t) = \sgn(c) [k^2 \sn (\varphi /k) (\cn (\varphi /k) - \cn (\varphi + t)/k) - \dn (\varphi /k) (t/k + \E(\varphi /k) - \E(\varphi + t)/k)].
\end{align*}

\i Geodesics of type Sep have the following expression :
\begin{align*}
&\cos \th(t) = 1/ (\cosh \varphi \cosh (\varphi + t))  + \tanh \varphi  \tanh (\varphi +t), \\
&\sin \th(t) = \sgn(\cos(\nu/2)) (\tanh \varphi /\cosh (\varphi +t) - \tanh (\varphi +t) /\cosh \varphi), \\
&x(t) = \sgn(\cos(\nu/2)) \sgn(c) [(1/\cosh \varphi)(1/\cosh \varphi - 1/\cosh (\varphi + t)) + \tanh \varphi(t + \tanh \varphi - \tanh (\varphi + t))],\\
&y(t) = \sgn(c) [\tanh \varphi (1/\cosh \varphi - 1/\cosh (\varphi + t)) -(1/\cosh \varphi) (t + \tanh \varphi - \tanh (\varphi + t))].
\end{align*}
\ei
Pictures of geodesics of type R, O, Sep are given in Figures \ref{fig:xyC1}, \ref{fig:xyC2} and \ref{fig:xyC3}, respectively.

\section{Proof of Theorem \ref{t-equiv}}
\label{a-pmp}

In this appendix, we prove Theorem \ref{t-equiv}. The structure of the proof is given in Section \ref{s-genpmp}. We are left to prove STEP 1 and STEP 2.\\

\noindent\b{STEP 1:} If $(\Gamma(.), v(.))$ is a solution of the generalized Pontryagin Maximum Principle for
\pcurve, then, the corresponding pair $(q(.),(u(.),v(.)))$ is a solution of the generalized Pontryagin Maximum
Principle for \pmec.

\bproof Without loss of generality, we provide the proof for $\xi=1$.

Apply the generalized PMP both to problems \pcurve\ and \pmec. For \pmec, the unmaximised Hamiltonian is $\mathcal{H}_M:=p_1\cos(\theta)u+p_2\sin(\theta)u+p_3v+\lam \sqrt{u^2+v^2}$. For \pcurve, replace $u$ with 1: we denote such Hamiltonian with $\mathcal{H}_C$.  We denote the maximised Hamiltonians with $H_M,H_C$, respectively. Recall that we study free time problems, thus both the maximised Hamiltonians satisfy $H_M\equiv 0$ and $H_C\equiv 0$, see \cite[Sec. 12.3]{agra-book}.

We observe that for both problems there are no strictly abnormal extremals (i.e. solutions with $\lam=0$). Indeed, for \pcurve\ abnormal extremals are straight lines, that can be realized as normal extremals too. The same holds for \pmec. Thus we fix from now on $\lam=-1$ without loss of generality.

Let now $(\bar \Gamma(.),\bar p(.),\bar v(.))$ be a trajectory vector-covector-control satisfying the generalized PMP for \pcurve. We prove that the corresponding trajectory vector-covector-controls $(\bar q(.),\bar p(.),(1,\bar v(.)))$ satisfies the generalized PMP for \pmec. The main point here is that $\mathcal{H}_M$ depends on two parameters $(u,v)$, while $\mathcal{H}_C$ depends on $v$ only. Thus, to maximise Hamiltonians, one has more degrees of freedom for $\mathcal{H}_M$ than for $\mathcal{H}_C$. We need to prove that such additional degree of freedom $u$ does not improve maximisation of the Hamiltonian.

We first prove that, if $\bar v(.)$ maximises\footnote{i.e., it maximises $\mathcal{H}_C$ along the trajectory $(\bar q(.),\bar p(.))$.}  $\mathcal{H}_C(\bar q(.),\bar p(.),v(.))$, then the choice $u(.)\equiv 1, v(.)=\bar v(.)$ maximises the Hamiltonian $\mathcal{H}_M(\bar q(.),\bar p(.),(u(.),v(.)))$. First observe that both $\mathcal{H}_M$ and $\mathcal{H}_C$ are $C^\infty$ (except for $\mathcal{H}_M$ in $(0,0)$), and concave with respect to variables $u,v$ and $v$, respectively. Moreover, we have no constraints on the controls. Thus, maximisation of the Hamiltonian is equivalent to have $\nabla_u \mathcal{H}=0$.

We are reduced to prove that $\frac{\partial \mathcal{H}_M}{\partial u}=\frac{\partial \mathcal{H}_M}{\partial v}=0$ when evaluated in $(\bar q(.),\bar p(.),(1,\bar v(.)))$. Observe that $\frac{\partial \mathcal{H}_M}{\partial v}=\frac{\partial \mathcal{H}_C}{\partial v}$ for $u=1$; thus, since $\bar v(.)$ maximises $\mathcal{H}_C$, then $\frac{\partial \mathcal{H}_C}{\partial v}=0$. Hence $\frac{\partial \mathcal{H}_M}{\partial v}=0$. A simple computation also shows that $ \frac{\partial \mathcal{H}_M}{\partial u}$ evaluated in 
 $u(.)\equiv 1, v(.)=\bar v(.)$ is $\bar p_1 \cos(\bar \theta)+\bar p_2 \sin(\bar \theta)-\frac{1}{\sqrt{1+\bar v^2}}$, whose expression coincides with $H_C$ when replacing $p_3$ with its expression with respect to the optimal control, that is $p_3=\frac{\bar v}{\sqrt{1+\bar v^2}}$. Since $H_C=0$, then $\frac{\partial \mathcal{H}_M}{\partial u}=0$, hence $\mathcal{H}_M$ is maximised by $u(.)\equiv 1, v(.)=\bar v(.)$.

Thus we have that $H_M=H_C$ on this trajectory. Then, since $H_C=0$, then it clearly holds $\mathcal{H}_M(\bar q(.),\bar p(.),(1,\bar v(.)))\equiv 0$ and it is also clear that $(\bar q(.),\bar p(.))$ is a solution of the Hamiltonian system with Hamiltonian $H_M$. Then, $(\bar q(.),\bar p(.))$ is a solution of the generalized PMP for \pmec.
\eproof\\

\noindent \b{STEP 2:} Let $(q(.),(u(.),v(.)))$ with $u(.)\equiv 1$ be a solution of the generalized Pontryagin Maximum
Principle for \pmec. Then, the curve reparametrized by \srarc\ is a solution of the standard Pontryagin Maximum Principle.

\bproof Recall that for \pmec\ one can always reparametrize curves by \srarc. This also transforms trajectoires with $L^1$ controls in trajectories with $L^\infty$ controls without changing the cost, as explained in Remark \ref{r-Linf}. Choose such reparametrization.

As a consequence, a solution to the generalized PMP can be reparametrized to have controls in $L^\infty$. Since the expression of the equations are the same for the standard and generalized PMP, then this reparametrized curve is a solution to the standard PMP.
\eproof

\newcommand{\auth}[1]{{\sc #1}}
\newcommand{\tit}[1]{{\rm #1}}
\newcommand{\titl}[1]{{\it #1}}
\newcommand{\jou}[1]{{\it #1}}
\newcommand{\vol}[1]{{\it #1}}
\newcommand{\pp}[1]{pp.~#1}

\end{document}